\documentclass[11pt]{article}
\usepackage{amssymb}
\usepackage{mathrsfs}
\usepackage{float}
\usepackage{amsmath}
\usepackage{amsfonts,amssymb}
\usepackage{graphicx,color}
\catcode`\@=11 \@addtoreset{equation}{section}

\catcode`\@=12

\newtheorem{Theorem}{Theorem}[section]
\newtheorem{Definition}{Definition}[section]

\newtheorem{Remark}{Remark}[section]

\def\2{{I \hskip -1.0mm I}}
\def\3{{I \hskip -1.0mm I\hskip -1.0mm I}}
\def\4{{I \hskip -0.9mm V}}
\def\6{{V \hskip -1.35mm I}}

\title{Dissipative Hyperbolic Geometric Flow}
\author{Wen-Rong Dai$^1$, \;\; De-Xing Kong$^1$
\;\; and \;\; Kefeng Liu$^{2,1}$\\
\\ $^{1}$Center of Mathematical
Sciences, Zhejiang University\\
 Hangzhou 310027, China\\
$^2$Department of Mathematics, UCLA, CA 90095, USA}
\date{ }
\begin{document}
\maketitle
\begin{abstract} In this paper we introduce a new kind of hyperbolic geometric flows ---
dissipative hyperbolic geometric flow. This kind of flow is defined
by a system of quasilinear wave equations with dissipative terms.
Some interesting exact solutions are given, in particular, a new
concept
--- hyperbolic Ricci soliton is introduced and some of its geometric
properties are described. We also establish the short-time existence
and uniqueness theorem for the dissipative hyperbolic geometric
flow, and prove the nonlinear stability of the flow defined on the
Euclidean space of dimension larger than 2. Wave character of the
evolving metrics and curvatures is illustrated and the nonlinear
wave equations satisfied by the curvatures are derived.

\vskip 6mm


\vskip 3mm

\end{abstract}
\newpage
\baselineskip=7mm

\section{Introduction}
Let $\mathscr{M}$  be an $n$-dimensional complete Riemannian
manifold with Riemannian metric $g_{ij}$. The following evolution
equation for the metric $g_{ij}$
\begin{equation}\frac{\partial^{2}g_{ij}}{\partial
t^{2}}+2R_{ij}+\mathscr{F}_{ij}\left(g,\frac{\partial g}{\partial
t}\right)=0\end{equation} has been recently introduced and named as
{\it general version of hyperbolic geometric flow} by Kong and Liu
\cite{kl}, where $R_{ij}$ is the corresponding Ricci curvature
tensor and $\mathscr{F}_{ij}$ is a given smooth symmetric tensor on
the Riemannian metric $g$ and its first order derivative with
respect to $t$. A special but important case is
\begin{equation}
\frac{\partial^{2}g_{ij}}{\partial t^{2}}=-2R_{ij}.
\end{equation}
Usually, we call (1.2) the {\it standard hyperbolic geometric flow}
or simply {\it hyperbolic geometric flow}. (1.1) and (1.2) are two
nonlinear systems of second order partial differential equations on
the metric $g_{ij}$.

For the hyperbolic geometric flow (1.2), some interesting exact
solutions have been constructed by Kong and Liu \cite{kl}. Recently,
Kong, Liu and Xu \cite{klx} have investigated the evolution of
Riemann surfaces under the flow (1.2) and given some results on the
global existence and blowup phenomenon of smooth solutions to the
flow equation (1.2). In our paper \cite{dkl}, we prove the
short-time existence for the hyperbolic geometric flow (1.2) and the
nonlinear stability of the Euclidean space with dimension larger
than 4. Moreover, we also study the wave character of the curvatures
for the flow (1.2) and derive the equations satisfied by curvatures
including the Riemannian curvature tensor $R_{ijkl}$, the Ricci
curvature tensor $R_{ij}$ and the scalar curvature $R$. However,
these evolution equations are quite complicated. In general, the
solution of the hyperbolic geometric flow (1.2) may blowup in a
finite time even for smooth initial data.

Motivated by the well-developed theory of the dissipative hyperbolic
equations, we introduce a new geometric analytical tool ---
dissipative hyperbolic geometric flow:
\begin{equation}
\begin{array}{lll}
\dfrac{\partial^{2}g_{ij}}{\partial t^{2}}& = &
-2R_{ij}+2g^{pq}\dfrac{\partial g_{ip}}{\partial t}\dfrac{\partial
g_{jq}}{\partial t}-\left(d+2g^{pq}\dfrac{\partial g_{pq}}{\partial
t}\right)\dfrac{\partial g_{ij}}{\partial t}+ \vspace{2mm}\\
& &\dfrac{1}{n-1}\left[\left(g^{pq}\dfrac{\partial g_{pq}}{\partial
t}\right)^2+ \dfrac{\partial g^{pq}}{\partial t}\dfrac{\partial
g_{pq}}{\partial t}\right]g_{ij}\end{array}
\end{equation}
where $g_{ij}(t)$ stands for a family of Riemannian metrics defined
on $\mathscr{M}$, and $d$ is a positive constant. The derivation of
(1.3) is given in Section 6. Here we would like to point out that
the reason that we choose (1.3) as the equation form of dissipative
hyperbolic geometric flow is that, in the case it possesses a
simpler equation satisfied by the scalar curvature. Noting the
dissipative property of (1.3), we expect that the dissipative
hyperbolic geometric flow admits a global smooth solution (i.e., a
family of Riemannian metrics) for all $t\ge 0$, and the solution
(metrics) has some good or anticipant geometric properties for
relatively general initial data in the case that the dissipative
coefficient $d$ is chosen to be suitably large.

In the present paper we will focus on some basic properties enjoyed
by the dissipative hyperbolic geometric flow. The first basic
property is on the {\it hyperbolic Ricci soliton}. The hyperbolic
Ricci soliton is a new concept which we introduce in this paper. We
will prove that there does not exist steady gradient hyperbolic
Ricci soliton with initial metric of positive average scalar
curvature on $n$-dimensional compact manifold (where $n\ge 3$).
Comparing with the traditional Ricci flow, here we need the
assumption that the initial metric has non-negative average scalar
curvature. If this assumption does not hold, then the question
whether there exist steady gradient hyperbolic Ricci solitons still
remains open. See Theorem 3.1 for the detail.

The second fundamental property is the short-time existence and
uniqueness theorem for the dissipative hyperbolic geometric flow.
For compact manifolds, we can prove that the dissipative hyperbolic
geometric flow always admits a unique smooth solution ( a family of
Riemannian metrics) for smooth initial data. See Theorem 4.1. Notice
that the dissipative hyperbolic geometric flow (1.3) is only weakly
hyperbolic, since the symbol of the derivative of
$E=E(g_{ij})\stackrel{\triangle}{=}-2R_{ij}$ has zero eigenvectors
in the natural coordinates. In order to reduce the nonlinear weakly
hyperbolic partial differential equation (1.3) to a nonlinear
symmetric system of strictly hyperbolic partial differential
equations, we use harmonic coordinates introduced by DeTurck and
Kazdan \cite{dk}. Then by the standard theory of symmetric
hyperbolic system, we can prove the short-time existence and
uniqueness theorem 4.1.

The third property is the nonlinear stability. By the global
existence theory of dissipative wave equations, we can prove the
global nonlinear stability of the Euclidean space $\mathbb{R}^n$
with $n\ge 3$. See Theorem 5.1 for the details. In the proof of
nonlinear stability, the dissipative property of the flow (1.3) play
an important role.

The fourth fundamental property is the wave character of the
curvatures. Since the dissipative hyperbolic geometric flow is
described by a system of quasilinear wave equations on the metrics
$g_{ij}(t,x)$, the wave property of the metric implies the wave
character of the curvatures. The equations will play an important
role in the future study. See Section 6 for the details.

The paper is organized as follows. In Section 2, we introduce the
dissipative hyperbolic geometric flow equation and give a useful
lemma. In order to understand the basics of the dissipative
hyperbolic geometric flow, we construct some exact solutions. These
solutions may be useful in physics. In Section 3, we introduce the
steady gradient hyperbolic Ricci soliton, and prove Theorem 3.1
--- one of the main results in this paper. Section 4 is devoted to the
short-time existence and uniqueness of the flow, while Section 5 is
devoted to the global nonlinear stability of the Euclidean space
$\mathbb{R}^n$ with $n\ge 3$. The wave character of the curvatures
is discussed in Section 6, and the nonlinear wave equations
satisfied by the curvatures are also derived in this section.

\section{Dissipative hyperbolic geometric flow}

The dissipative hyperbolic geometric flow considered here is defined
by the equation (1.3), namely,
\begin{equation}
\begin{array}{lll}
\dfrac{\partial^{2}g_{ij}}{\partial t^{2}}& = &
-2R_{ij}+2g^{pq}\dfrac{\partial g_{ip}}{\partial t}\dfrac{\partial
g_{jq}}{\partial t}-\left(d+2g^{pq}\dfrac{\partial g_{pq}}{\partial
t}\right)\dfrac{\partial g_{ij}}{\partial t}+ \vspace{2mm}\\
& &\dfrac{1}{n-1}\left[\left(g^{pq}\dfrac{\partial g_{pq}}{\partial
t}\right)^2+ \dfrac{\partial g^{pq}}{\partial t}\dfrac{\partial
g_{pq}}{\partial t}\right]g_{ij}\end{array}
\end{equation}
where $g_{ij}(t)$ stands for a family of Riemannian metrics defined
on $\mathscr{M}$, and $d$ is a positive constant. The reason that we
choose (2.1) as the equation form of dissipative hyperbolic
geometric flow is as follows: in this case the flow possesses a
simpler equation satisfied by the scalar curvature. See the
derivation of (2.1) in Section 6.

We first establish some useful equations from the flow equation
(2.1). Let
\begin{eqnarray}
u(x,t)&=&g^{ij}\frac{\partial g_{ij}}{\partial t},\vspace{2mm}\\
v(x,t)&=&\left|\frac{\partial g}{\partial
t}\right|^2=g^{ik}g^{jl}\frac{\partial g_{ij}}{\partial
t}\frac{\partial g_{kl}}{\partial
t},\vspace{2mm}\\w(x,t)&=&g^{ik}g^{jl}g^{pq}\frac{\partial
g_{ip}}{\partial t}\frac{\partial g_{jq}}{\partial t}\frac{\partial
g_{kl}}{\partial t}
\end{eqnarray}
and denote the matrix
\begin{equation}
G(x,t)=\left(\frac{\partial g_{ij}}{\partial
t}g^{jk}\right).\end{equation} Then we have
\begin{eqnarray}
u(x,t)=trG(x,t),\;\;v(x,t)=trG^2(x,t),\;\;w(x,t)=trG^3(x,t),
\end{eqnarray}
where $trG$ stands for the trace of the matrix $G$. Thus by (2.1) we
obtain
\begin{equation}\begin{array}{lll} {\displaystyle \frac{\partial u(x,t)}{\partial
t}}&=&{\displaystyle \frac{\partial }{\partial
t}\left(g^{ij}\frac{\partial g_{ij}}{\partial
t}\right)}\vspace{2mm}\\ &=&{\displaystyle \frac{\partial
g^{ij}}{\partial t}\frac{\partial g_{ij}}{\partial
t}+g^{ij}\frac{\partial^2
g_{ij}}{\partial t^2}}\vspace{2mm}\\
&=&{\displaystyle -g^{ik}g^{jl}\frac{\partial g_{ij}}{\partial
t}\frac{\partial g_{kl}}{\partial
t}+g^{ij}\left[-2R_{ij}+2g^{pq}\frac{\partial g_{ip}}{\partial
t}\frac{\partial g_{jq}}{\partial t}-2g^{pq}\frac{\partial
g_{pq}}{\partial t}\frac{\partial
g_{ij}}{\partial t}\right.}\vspace{2mm}\\
& &{\displaystyle \left. -d\cdot\frac{\partial g_{ij}}{\partial
t}+\frac{1}{n-1}\left(\left(g^{pq}\frac{\partial g_{pq}}{\partial
t}\right)^2+\left(\frac{\partial g^{pq}}{\partial t}\frac{\partial
g_{pq}}{\partial
t}\right)\right)g_{ij}\right]}\vspace{2mm}\\&=&{\displaystyle
-2R-\frac{n-2}{n-1}u^2-du-\frac{1}{n-1}v}
\end{array}\end{equation}
and
\begin{equation}\begin{array}{lll} {\displaystyle \frac{\partial v(x,t)}{\partial t}}&=& {\displaystyle 2\frac{\partial
g^{ij}}{\partial t}g^{pq}\frac{\partial g_{ip}}{\partial
t}\frac{\partial g_{jq}}{\partial t}+2g^{ij}g^{pq}\frac{\partial^2
g_{ip}}{\partial t^2}\frac{\partial g_{jq}}{\partial
t}}\vspace{2mm}\\ &=&{\displaystyle
-2g^{ir}g^{js}g^{pq}\frac{\partial g_{ip}}{\partial t}\frac{\partial
g_{jq}}{\partial t}\frac{\partial g_{rs}}{\partial
t}+2g^{ij}g^{pq}\frac{\partial g_{jq}}{\partial
t}\left[-2R_{ip}+2g^{rs}\frac{\partial g_{ir}}{\partial
t}\frac{\partial g_{ps}}{\partial t}\right.}\vspace{2mm}\\
&& {\displaystyle \left.-2\left(g^{rs}\frac{\partial
g_{rs}}{\partial t}\right)\frac{\partial g_{ip}}{\partial
t}-d\frac{\partial g_{ip}}{\partial
t}+\frac{1}{n-1}\left(g^{rs}\frac{\partial g_{rs}}{\partial
t}\right)^2g_{ip}+\frac{1}{n-1}\left(\frac{\partial g^{rs}}{\partial
t}\frac{\partial g_{rs}}{\partial
t}\right)g_{ip}\right]}\vspace{2mm}\\&=& {\displaystyle
2w-4g^{ik}g^{jl}\frac{\partial g_{ij}}{\partial
t}R_{kl}-\left(4+\frac{2}{n-1}\right)uv-2dv+\frac{2}{n-1}u^3.}\end{array}\end{equation}

\begin{Theorem}
For the dissipative hyperbolic geometric flow (2.1), the quantities
$u(x,t)$, $v(x,t)$ and $w(x,t)$ satisfy the following equations
\begin{eqnarray}
\frac{\partial u(x,t)}{\partial
t}=-2R-\frac{n-2}{n-1}u^2-du-\frac{1}{n-1}v
\end{eqnarray}
and
\begin{eqnarray}
\frac{\partial v(x,t)}{\partial t}=2w-4g^{ik}g^{jl}\frac{\partial
g^{ij}}{\partial t}R_{kl}-(4+\frac{2}{n-1})uv-2dv+\frac{2}{n-1}u^3.
\end{eqnarray}
\end{Theorem}

In order to understand basically the dissipative hyperbolic
geometric flow, in what follows we construct some exact solutions.

Consider the following Cauchy problem
\begin{equation}\left\{\!\!\!\!\!\!
\begin{array}{lll}
&\dfrac{\partial^{2}g_{ij}}{\partial t^{2}} = &
-2R_{ij}+2g^{pq}\dfrac{\partial g_{ip}}{\partial t}\dfrac{\partial
g_{jq}}{\partial t}-\left(d+2g^{pq}\dfrac{\partial g_{pq}}{\partial
t}\right)\dfrac{\partial g_{ij}}{\partial t} \vspace{2mm}\\
& & + \dfrac{1}{n-1}\left[\left(g^{pq}\dfrac{\partial
g_{pq}}{\partial t}\right)^2+ \dfrac{\partial g^{pq}}{\partial
t}\dfrac{\partial
g_{pq}}{\partial t}\right]g_{ij}, \vspace{2mm}\\
& g_{ij}(x,0) & =g^0_{ij}(x), {\displaystyle \quad \frac{\partial
g_{ij}}{\partial t}(x,0)=k^0_{ij}(x),}
\end{array}\right.
\end{equation}
where $g_{ij}^0(x)$ is a Riemannian metric on the manifold
$\mathscr{M}$, and $k^0_{ij}(x)$ is a symmetric tensor on
$\mathscr{M}$.

If we assume that the initial metric $g^0_{ij}(x)$ is Ricci flat,
and the initial velocity $k^0_{ij}(x)$ vanishes, then easily see
that $g_{ij}(x,t)= g^0_{ij}(x)$ is the unique smooth solution to the
Cauchy problem (2.11).

If we assume that the initial Riemannian metric is Einstein, that is
to say,
\begin{eqnarray}
R_{ij}(x,0)=\lambda g_{ij}(x,0),\;\;\forall\;x\in\mathscr{M},
\end{eqnarray}
where $\lambda$ is a constant. Furthermore, we suppose that
\begin{eqnarray}
\frac{\partial g_{ij}}{\partial t}(x,0)=\mu g_{ij}(x,0),
\end{eqnarray}
where $\mu$ is an another constant. Let
\begin{eqnarray}
g_{ij}(x,t)=\rho(t)g_{ij}(x,0).
\end{eqnarray}
By the definition of the Ricci tensor, we have
\begin{eqnarray}
R_{ij}(x,t)=R_{ij}(x,0)=\lambda
g_{ij}(x,0),\;\;\forall\;x\in\mathscr{M}.
\end{eqnarray}
It follows from (2.13) and (2.14) that
\begin{eqnarray}
\rho(0)=1,\quad\rho^{\prime}(0)=\mu.
\end{eqnarray}
Substituting (2.14) into the evolution equation (2.1) gives the
following ODE
\begin{eqnarray}
\rho^{\prime\prime}(t)=-d\rho^{\prime}(t)-2\lambda.
\end{eqnarray}
The solution of (2.17) with the initial data (2.16) reads
\begin{eqnarray}
\rho(t)=1-\frac{2\lambda}{d}t-\left(\frac{\mu}{d}+\frac{2\lambda}{d^2}\right)(e^{-dt}-1).
\end{eqnarray}
It follows from (2.18) that
\begin{eqnarray}
\rho^{\prime}(t)=-\frac{2\lambda}{d}+\left(\mu+\frac{2\lambda}{d}\right)e^{-dt}.
\end{eqnarray}
Noting that $d>0$, we distinguish the following three cases to
discuss:

{\bf Case I.} $\;\;\lambda>0.$

In this case, it follows from (2.18) that
$$\lim_{t\rightarrow
+\infty}\rho(t)=-\infty.$$ Thus the evolving metric $g_{ij}(x,t)$
shrinks homothetically to a point as $t$ approaches some finite time
$T$.

{\bf Case II.} $\;\;\lambda=0.$

In the present situation, $\rho(t)=1-\frac{\mu}{d}(e^{-dt}-1)$. If
$\frac{\mu}{d}<-1$, then the evolving metric $g_{ij}(x,t)$ shrinks
homothetically to a point as $t$ approaches the time
$T\stackrel{\triangle}{=} -\frac{1}{d}\ln(1+\frac{d}{\mu})$; If
$\frac{\mu}{d}>-1$, then the metric $g_{ij}(x,t)$ evolves smoothly
and is positive defined for all time; If $\frac{\mu}{d}=-1$, the
metric $g_{ij}(x,t)$ evolves smoothly and is positive defined for
all time, but it shrinks homothetically to a point as $t\rightarrow
+\infty$.

{\bf Case III.} $\;\;\lambda < 0.$

In this case, if $\mu<0$ and $\rho(T_0)\le 0,$ where
$T_0\stackrel{\triangle}{=}-\frac{1}{d}\ln\left(\frac{2\lambda}{2\lambda+d\mu}\right)$,
then the evolving metric $g_{ij}(x,t)$ shrinks homothetically to a
point as $t$ approaches some finite time not later than $T$.
Otherwise, $g_{ij}(x,t)$ is smooth and positive defined for all
time.

Summarizing the above argument leads to the following theorem.

\begin{Theorem}
For the Cauchy problem (2.11) of the dissipative hyperbolic
geometric flow, suppose that the assumptions (2.12)-(2.13) are
satisfied. Then, if one of the following conditions is satisfied,
then the evolving metric $g_{ij}(x,t)$ shrinks homothetically to a
point as $t$ approaches some finite time:

(a) $\;\;\lambda>0$;

(b) $\;\;\lambda=0$ and $\mu<-d$;

(c) $\;\;\lambda<0$, $\mu<0$ and
$\rho\left(\frac{1}{d}\ln\left(\frac{2\lambda}{2\lambda+d\mu}\right)\right)\ge
0$.

For the other instances, $g_{ij}(x,t)$ are smooth and positive
defined for all time. In addition, if $\lambda=0$ and $\mu=-d<0$,
the metric $g_{ij}(x,t)$ evolves smoothly and is positively defined
for all time, but it shrinks homothetically to a point as
$t\rightarrow +\infty$.
\end{Theorem}

\section{Hyperbolic Ricci soliton}

The theory of soliton solutions plays an important role in the study
of geometric analysis, in particular in the study of Ricci flow. In
this section we first introduce a new concept --- steady hyperbolic
Ricci soliton for the flow (2.1), and then describe its properties.

\begin{Definition} A solution to an evolution
equation is called a steady soliton, if it evolves under a
one-parameter subgroup of the symmetry group of the equation; A
solution to the dissipative hyperbolic geometric flow (2.1)  is
called a steady hyperbolic Ricci soliton, if it moves by a
one-parameter subgroup of the symmetry group of the equation (2.1).
\end{Definition}

If $\varphi_t$ is a one-parameter group of diffeomorphisms generated
by a vector field $V$ on $\mathscr{M}$, then the hyperbolic Ricci
soliton is given by
\begin{equation}
g_{ij}(x,t)=\varphi_t^*g_{ij}(x,0)=g_{ij}(\varphi_t(x),0).
\end{equation}
It implies that
\begin{equation}
\frac{\partial}{\partial
t}g_{ij}(x,t)=\mathfrak{L}_Vg_{ij}=g_{ik}\nabla_jV^k+g_{jk}\nabla_iV^k\triangleq
T_{ij}
\end{equation}
and \begin{equation}\begin{array}{lll} \frac{\partial}{\partial
t^2}g_{ij}(x,t)&=&\mathfrak{L}_V\mathfrak{L}_Vg_{ij}=\mathfrak{L}_VT_{ij}\vspace{2mm}\\
&=&T_{ij;k}V^k+T_{kj}V^k_{;i}+T_{ki}V^k_{;j}\vspace{2mm}\\
&=&(g_{ip}\nabla_jV^p+g_{jp}\nabla_iV^p)_{;k}V^k+(g_{kp}\nabla_jV^p+g_{jp}\nabla_kV^p)V^k_{;i}
\vspace{2mm}\\
& & +(g_{kp}\nabla_iV^p+g_{ip}\nabla_kV^p)V^k_{;j}\vspace{2mm}\\
&=&\left(g_{ip}\nabla_k\nabla_jV^p+g_{jp}\nabla_k\nabla_iV^p\right)V^k
+g_{kp}(\nabla_iV^k\cdot\nabla_jV^p+\nabla_jV^k\cdot\nabla_iV^p)\vspace{2mm}\\
&
&+g_{ip}\nabla_jV^k\cdot\nabla_kV^p+g_{jp}\nabla_iV^k\cdot\nabla_kV^p,
\end{array}\end{equation}
where $\mathfrak{L}_V$ stands for the Lie derivative with respect to
the vector field $V$. Thus, the equation (2.1) can be reduced to
\begin{eqnarray}
\nonumber&&\left(g_{ip}\nabla_k\nabla_jV^p+g_{jp}\nabla_k\nabla_iV^p\right)V^k
+g_{kp}(\nabla_iV^k\cdot\nabla_jV^p+\nabla_jV^k\cdot\nabla_iV^p)\vspace{2mm}\\
\nonumber&&+g_{ip}\nabla_jV^k\cdot\nabla_kV^p+g_{jp}\nabla_iV^k\cdot\nabla_kV^p
\vspace{2mm}\\
\nonumber&=&-2R_{ij}+2g^{pq}(g_{ik}\nabla_pV^k+g_{pk}\nabla_iV^k)(g_{jl}\nabla_qV^l+g_{ql}\nabla_jV^l)\vspace{2mm}\\
\nonumber&&-2g^{pq}(g_{pk}\nabla_qV^k+g_{qk}\nabla_pV^k)(g_{il}\nabla_jV^l+g_{jl}\nabla_iV^l)
-d(g_{ik}\nabla_jV^k+g_{jk}\nabla_iV^k)\vspace{2mm}\\
\nonumber&&+\frac{1}{n-1}\left[g^{pq}(g_{pk}\nabla_qV^k+g_{qk}\nabla_pV^k)\right]^2g_{ij}\vspace{2mm}\\
&&-\frac{1}{n-1}\left[g^{pr}g^{qs}(g_{pk}\nabla_qV^k+g_{qk}\nabla_pV^k)
(g_{rl}\nabla_sV^l+g_{sl}\nabla_rV^l)\right]g_{ij}.
\end{eqnarray}
We predigest it into the following
\begin{eqnarray}
\nonumber&&2R_{ij}+\left(g_{ip}\nabla_k\nabla_jV^p+g_{jp}\nabla_k\nabla_iV^p\right)V^k\vspace{2mm}\\
\nonumber&=&2g^{pq}g_{ik}g_{jl}\nabla_pV^k\nabla_qV^l+g_{ik}\nabla_jV^l\nabla_lV^k
+g_{jk}\nabla_iV^l\nabla_lV^k\vspace{2mm}\\
\nonumber&&- (d+4\nabla_kV^k)(g_{il}\nabla_jV^l+g_{jl}\nabla_iV^l)+\frac{4}{n-1}(\nabla_qV^q)^2g_{ij}\vspace{2mm}\\
&&-\frac{2}{n-1}\left(g_{kl}g^{pq}\nabla_pV^k\nabla_qV^l+\nabla_pV^q\nabla_qV^p\right)g_{ij}.
\end{eqnarray}

If the vector field $V$ is the gradient of a function $f$ on
$\mathscr{M}$, then the soliton is called {\bf a steady gradient
hyperbolic Ricci soliton}. In what follows, we consider the steady
gradient hyperbolic Ricci soliton.

For the steady gradient hyperbolic Ricci soliton, the equation (3.5)
becomes
\begin{eqnarray*}
&&2R_{ij}+\left(g_{ip}\nabla_k\nabla_j\nabla^pf+g_{jp}\nabla_k\nabla_i\nabla^pf\right)\nabla^kf\vspace{2mm}\\
&=&2g^{pq}g_{ik}g_{jl}\nabla_p\nabla^kf\nabla_q\nabla^lf+g_{ik}\nabla_j\nabla^lf\nabla_l\nabla^kf
+g_{jk}\nabla_i\nabla^lf\nabla_l\nabla^kf\vspace{2mm}\\
&& -(d+4\nabla_k\nabla^kf)(g_{il}\nabla_j\nabla^lf+g_{jl}\nabla_i\nabla^lf)+\frac{4}{n-1}(\nabla_q\nabla^qf)^2g_{ij}\vspace{2mm}\\
&&-\frac{2}{n-1}\left(g_{kl}g^{pq}\nabla_p\nabla^kf\nabla_q\nabla^lf+\nabla_p\nabla^qf\nabla_q\nabla^pf\right)g_{ij}.
\end{eqnarray*}
That is to say,
\begin{eqnarray}
\nonumber
R_{ij}+\nabla_k(\nabla_i\nabla_jf)\nabla^kf&=&2g^{pq}\nabla_p\nabla_if\nabla_q\nabla_jf-(d+4\triangle
f)\nabla_i\nabla_jf
\vspace{2mm}\\
&&+\frac{2}{n-1}(\triangle
f)^2g_{ij}-\frac{2}{n-1}(g^{pq}g^{kl}\nabla_p\nabla_kf\nabla_q\nabla_lf)g_{ij}.
\end{eqnarray}
Taking the trace on $i$ and $j$ yields
\begin{equation}
R+\nabla_k(\triangle f\cdot \nabla^k
f)=-\frac{2}{n-1}|\nabla^2f|^2-\frac{n-3}{n-1}(\triangle
f)^2-d\cdot\triangle f.
\end{equation}
Thus, the following theorem comes easily from (3.5)-(3.7).
\begin{Theorem}
For the dissipative hyperbolic geometric flow,  (3.5) and (3.6) are
the evolution equations satisfied by the steady hyperbolic Ricci
soliton and the steady gradient hyperbolic Ricci soliton,
respectively. Furthermore, for an $n$-dimensional compact manifold
with $n\ge 3$, if the average scalar curvature of the initial metric
is non-negative, i.e.,
\begin{eqnarray}
r(0)\triangleq \dfrac
{\int_{\mathscr{M}}R(x,0)dV}{\int_{\mathscr{M}}dV}\ge 0,
\end{eqnarray}
then for the steady gradient hyperbolic Ricci soliton, the
generating function $f$ must satisfy the condition $Hess(f)\equiv 0
\quad \mbox {on} \;\;\mathscr{M}$, i.e., $f$ is a constant and the
solution metric $g_{ij}(x,t)\equiv g_{ij}(x,0)$ is Ricci flat for
all time $t$. In reverse, if the initial metric $g_{ij}(x,0)$ is
Ricci flat and the function $f\equiv constant$, then it is obvious
that the steady gradient hyperbolic Ricci soliton generated by $f$
is a solution to the dissipative hyperbolic geometric flow.
\end{Theorem}

\section{Short-time existence and uniqueness}

In this section, we reduce the dissipative hyperbolic geometric flow
(2.1) to a symmetric hyperbolic system in the so-called harmonic
coordinates (see \cite{dk}), then based on this, we prove the
short-time existence and uniqueness theorem for the flow equation
(2.1).

Let $g_{ij}(x,t)$ be a family of metrics on an $n>1$ dimensional
manifold $\mathscr{M}$. We consider the space-time
$\mathbb{R}\times\mathscr{M}$ equipped with the following Lorentzian
metric
\begin{eqnarray}
ds^2=-dt^2+g_{ij}(x,t)dx^idx^j.
\end{eqnarray}
It follows from (3.4) in Dai, Kong and Liu \cite{dkl} that
\begin{eqnarray}
\nonumber\dfrac{\partial^{2} g_{ij}}{\partial
t^{2}}+2R_{ij}&=&\dfrac{\partial ^{2}g_{ij}}{\partial
t^{2}}-g^{kl}\dfrac{\partial^{2}g_{ij}}{\partial x^{k}\partial
x^{l}}+\left(g_{ik}\dfrac{\partial \Gamma^{k}}{\partial
x^{j}}+g_{jk}\dfrac{\partial \Gamma^{k}}{\partial x^{i}}\right)\vspace{2mm}\\
\nonumber & &+2g^{kl}g_{pq}\Gamma
_{ik}^{p}\Gamma_{jl}^{q}+\dfrac{\partial g_{ij}}{\partial
x_k}\Gamma^k\vspace{2mm}\\ &&+\left(g_{ik}\Gamma
_{rs}^{k}g^{pr}g^{qs}\dfrac{\partial g_{pq}}{\partial
x^{j}}+g_{jk}\Gamma _{rs}^{k}g^{pr}g^{qs}\dfrac{\partial
g_{pq}}{\partial x^i}\right),
\end{eqnarray}
where
\begin{eqnarray}
\Gamma^k\triangleq g^{ij}\Gamma_{ij}^k.
\end{eqnarray}
Then the evolution equation (2.1) for the dissipative hyperbolic
geometric flow can be reduced to the following
\begin{eqnarray}
\nonumber\dfrac{\partial^{2} g_{ij}}{\partial
t^{2}}-g^{kl}\dfrac{\partial^{2}g_{ij}}{\partial x^{k}\partial
x^{l}}&=&-\left(g_{ik}\dfrac{\partial \Gamma^{k}}{\partial
x^{j}}+g_{jk}\dfrac{\partial \Gamma^{k}}{\partial
x^{i}}\right)-2g^{kl}g_{pq}\Gamma
_{ik}^{p}\Gamma_{jl}^{q}-\dfrac{\partial g_{ij}}{\partial
x_k}\Gamma^k \vspace{2mm}\\ \nonumber&&-\left(g_{ik}\Gamma
_{rs}^{k}g^{pr}g^{qs}\dfrac{\partial g_{pq}}{\partial
x^{j}}+g_{jk}\Gamma _{rs}^{k}g^{pr}g^{qs}\dfrac{\partial
g_{pq}}{\partial x^i}\right)\vspace{2mm}\\
\nonumber&&+2g^{pq}\frac{\partial g_{ip}}{\partial t}\frac{\partial
g_{jq}}{\partial t}-2g^{pq}\frac{\partial g_{pq}}{\partial
t}\frac{\partial g_{ij}}{\partial t}-d\frac{\partial
g_{ij}}{\partial t}\vspace{2mm}\\& &
+\frac{1}{n-1}(g^{pq}\frac{\partial g_{pq}}{\partial
t})^2g_{ij}+\frac{1}{n-1}(\frac{\partial g^{pq}}{\partial
t}\frac{\partial g_{pq}}{\partial t})g_{ij}.
\end{eqnarray}

Similar to \cite{dk}, we make use of the harmonic coordinates such
that, for fixed time $t$, it holds that
\begin{eqnarray}
\Gamma^k(x,t)\triangleq g^{ij}\Gamma_{ij}^k\equiv 0, \;\;\mbox{when
$x$ is in an open neighborhood of point} \;\;p\in \mathscr{M}.
\end{eqnarray}
Then the equation (4.4) can be written as
\begin{eqnarray}
\dfrac{\partial^{2} g_{ij}}{\partial
t^{2}}=g^{kl}\dfrac{\partial^{2}g_{ij}}{\partial x^{k}\partial
x^{l}}+\widetilde{H_{ij}}(g_{kl},\frac{\partial g_{kl}}{\partial
t},\frac{\partial g_{kl}}{\partial x^p}),
\end{eqnarray}
where
\begin{eqnarray}
\nonumber\widetilde{H_{ij}}(g_{kl},\frac{\partial g_{kl}}{\partial
t},\frac{\partial g_{kl}}{\partial x^p})&=&-2g^{kl}g_{pq}\Gamma
_{ik}^{p}\Gamma_{jl}^{q}-\left(g_{ik}\Gamma
_{rs}^{k}g^{pr}g^{qs}\dfrac{\partial g_{pq}}{\partial
x^{j}}+g_{jk}\Gamma _{rs}^{k}g^{pr}g^{qs}\dfrac{\partial
g_{pq}}{\partial x^i}\right)\vspace{2mm}\\
\nonumber&&+2g^{pq}\frac{\partial g_{ip}}{\partial t}\frac{\partial
g_{jq}}{\partial t}-2g^{pq}\frac{\partial g_{pq}}{\partial
t}\frac{\partial g_{ij}}{\partial t}-d\frac{\partial
g_{ij}}{\partial t}\vspace{2mm}\\& &
+\frac{1}{n-1}(g^{pq}\frac{\partial g_{pq}}{\partial
t})^2g_{ij}+\frac{1}{n-1}(\frac{\partial g^{pq}}{\partial
t}\frac{\partial g_{pq}}{\partial t})g_{ij}
\end{eqnarray}
are homogenous quadratic with respect to $\frac{\partial
g_{kl}}{\partial x^p}$ and $\frac{\partial g_{kl}}{\partial t}$
except the dissipative term $d\frac{\partial g_{ij}}{\partial t}$
and rational with respect to $g_{kl}$ with non-zero denominator
$det(g_{ij})\neq 0$. By introducing the new unknowns
$g_{ij},\;\;h_{ij}=\frac{\partial g_{ij}}{\partial
t},\;\;g_{ij,k}=\frac{\partial g_{kl}}{\partial x^k}$, the system
(4.6) can be transformed into a system of partial differential
equations of first order
\begin{eqnarray}
\begin{cases}
\dfrac{\partial g_{ij}}{\partial t}=h_{ij}, \vspace{2mm}\\
g^{kl}\dfrac{\partial g_{ij,k}}{\partial t}=g^{kl}\dfrac{\partial
h_{ij}}{\partial x^{k}},\vspace{2mm}\\
\dfrac{\partial h_{ij}}{\partial t}=g^{kl}\dfrac{\partial
g_{ij,k}}{\partial x^{l}}+\widetilde{H}_{ij}.
\end{cases}
\end{eqnarray}
In the $C^2$ class, the system (4.8) is equivalent to (4.6). It is
easy to see that (4.8) is a quasilinear symmetric hyperbolic system,
which can be rewritten as
\begin{equation}
A^0(u)\dfrac{\partial u}{\partial t}=A^j(u)\dfrac{\partial
u}{\partial x^j}+B(u),
\end{equation}
where $u=(g_{ij},g_{ij,k},h_{ij})^T$ is the
$\dfrac{1}{2}n(n+1)(n+2)$-dimensional unknown vector function and
the coefficient matrices $A^0,A^j,B$ are given by
$$
A^0(u)=A^0(g_{ij},g_{ij,k},h_{ij})=\left(\begin{array}{cccccc}I & 0 & 0 & \cdots & 0 & 0\\
 0 & g^{11}I & g^{12}I & \cdots & g^{1n}I & 0\\
 0 & g^{21}I & g^{22}I & \cdots & g^{2n}I & 0\\
 \vdots & \cdots\\\texttt{}
 0 & g^{n1}I & g^{n2}I & \cdots & g^{nn}I & 0\\
 0 & 0 & 0 & \cdots & 0 & I\end{array}\right),
$$
$$
A^{j}(u)=A^{j}(g_{kl},g_{kl,p},h_{kl})=\left(\begin{array}{cccccc}0
& 0 & 0 & \cdots & 0 & 0\\
0 & 0 & 0 & \cdots & 0 & g^{j1}I\\
0 & 0 & 0 & \cdots & 0 & g^{j2}I\\
\cdots & \cdots &\cdots\\
0 & 0 & 0 & \cdots & 0 & g^{jn}I\\
0 & g^{1j}I & g^{2j}I & \cdots & g^{nj}I & 0\end{array}\right),
$$
where 0 is the
$\left(\dfrac{1}{2}n(n+1)\right)\times\left(\dfrac{1}{2}n(n+1)\right)$
zero matrix, $I$ is the
$\left(\dfrac{1}{2}n(n+1)\right)\times\left(\dfrac{1}{2}n(n+1)\right)$
 identity matrix,
$$B(u)=B(g_{ij},g_{ij,p},h_{ij})=\left(\begin{array}{c}h_{ij}\\
0\\
 \widetilde{H}_{ij}\\\end{array}\right),$$
in which 0 is the $\dfrac{1}{2}n^2(n+1)$-dimensional zero vector.

By the theory of the symmetric hyperbolic system ( \cite{fm},
\cite{fr}), we can obtain the following theorem.
\begin{Theorem}
Let $(\mathscr{M},\;g_{ij}^0(x))$ be an $n$-dimensional compact
Riemannian manifold. Then there exists a constant $\eta>0$ such that
the Cauchy problem (2.11) has a unique smooth solution $g_{ij}(x,t)$
on $\mathscr{M}\times [0,\eta]$.\end{Theorem}

\begin{Remark} Theorem 4.1 can also be proved in a manner similar to that in
DeTurck (see \cite{dt}, \cite{cz}, \cite{dkl}).
\end{Remark}

\section{Nonlinear stability of Euclidean metrics}

This section is devoted to the nonlinear stability of the
dissipative hyperbolic geometric flow (2.1) defined on the Euclidean
space with the dimension larger than two.

We consider the following Cauchy problem for the dissipative
hyperbolic geometric flow (2.1),
\begin{equation}\left\{\!\!\!\!\!\!
\begin{array}{lll}
&\dfrac{\partial^{2}g_{ij}}{\partial t^{2}} & =
-2R_{ij}+2g^{pq}\dfrac{\partial g_{ip}}{\partial t}\dfrac{\partial
g_{jq}}{\partial t}-\left(d+2g^{pq}\dfrac{\partial g_{pq}}{\partial
t}\right)\dfrac{\partial g_{ij}}{\partial t} \vspace{2mm}\\
& & \;\;\;\;+ \dfrac{1}{n-1}\left[\left(g^{pq}\dfrac{\partial
g_{pq}}{\partial t}\right)^2+ \dfrac{\partial g^{pq}}{\partial
t}\dfrac{\partial
g_{pq}}{\partial t}\right]g_{ij}, \vspace{2mm}\\
& g_{ij}(x,0) & {\displaystyle =\delta_{ij}+\epsilon g_{ij}^0(x),\;
\frac{\partial g_{ij}}{\partial t}(x,0)=\epsilon g_{ij}^1(x),}
\end{array}\right.
\end{equation}
where $g_{ij}^0(x)$ and $g_{ij}^1(x)$ are given symmetric tensors
defined on the Euclidean space $\mathbb{R}^n$.
\begin{Theorem}
The flat metric $g_{ij}=\delta_{ij}$ on the Euclidean space
$\mathbb{R}^n$ with $n\ge 3$ is globally nonlinearly stable with
respect to the given tensor $\left(g_{ij}^0(x),g_{ij}^1(x)\right)\in
C_0^{\infty}(\mathbb{R}^n)$ , i.e., there exists a positive constant
$\epsilon_0=\epsilon_0\left(g_{ij}^0(x),g_{ij}^1(x)\right)>0$ such
that, for any $\epsilon\in (0,\epsilon_0]$, the initial value
problem (4.1) admits a unique smooth solution $g_{ij}(x,t)$ for all
time $t\ge 0$.
\end{Theorem}
\begin{Remark}
For the standard hyperbolic geometric flow (1.2), we can only obtain
the nonlinear stability of the Euclidean space $\mathbb{R}^n$ with
$n\ge 5$ (see \cite{dkl}). Under suitable assumptions, similar
results are true for general hyperbolic geometric flow (1.1).
\end{Remark}

\noindent\textbf{Proof of Theorem 5.1.}  Let the symmetric tensor
$h_{ij}$ on $\mathbb{R}^n$ defined by
\begin{eqnarray}
h_{ij}(x,t)=g_{ij}(x,t)-\delta_{ij}
\end{eqnarray}
and $\delta^{ij}$ be the inverse of $\delta_{ij}$. Then for small
$h$,
\begin{eqnarray}
H^{ij}\triangleq g^{ij}-\delta^{ij}=-h^{ij}+O^{ij}(h^2),
\end{eqnarray}
where $h^{ij}=\delta^{ik}\delta^{jl}h_{kl}$ and $O^{ij}(h^2)$
vanishes to the second order at $h=0$. Then the Cauchy problem (5.1)
for the metric $g_{ij}(x,t)$ is equivalent to the following initial
 value problem for the tensor $h_{ij}(x,t)$ in the harmonic coordinates
${x^i}$ around the origin in $\mathbb{R}^n$
\begin{equation}
\begin{cases}{\displaystyle
\frac{\partial^2}{\partial
t^2}h_{ij}(x,t)=(\delta^{kl}+H^{kl})\frac{\partial^2h_{ij}}{\partial
x^kx^l}+\widetilde{H_{ij}}(\delta_{kl}+h_{kl},\frac{\partial
h_{kl}}{\partial t},\frac{\partial h_{kl}}{\partial x^p}),}
\\{\displaystyle t=0:h_{ij}(x,0)=\epsilon g_{ij}^0(x),\ \frac{\partial
h_{ij}}{\partial t}(x,0)=\epsilon g_{ij}^1(x),}
\end{cases}\end{equation}
where $\widetilde{H_{ij}}(\delta_{kl}+h_{kl},\frac{\partial
h_{kl}}{\partial t},\frac{\partial h_{kl}}{\partial x^p})$ is
defined in (4.7). Thus, the Cauchy problem (5.4) can be reduced to
the following
\begin{equation}
\begin{cases}{\displaystyle
\frac{\partial^2}{\partial
t^2}h_{ij}(x,t)-\delta^{kl}\frac{\partial^2h_{ij}}{\partial
x^kx^l}+d\frac{\partial h_{ij}}{\partial
t}=\bar{H}_{ij}(\delta_{kl}+h_{kl},\frac{\partial h_{kl}}{\partial
t},\frac{\partial h_{kl}}{\partial x^p}),}
\\{\displaystyle t=0:h_{ij}(x,0)=\epsilon g_{ij}^0(x),\ \frac{\partial
g_{ij}}{\partial t}(x,0)=\epsilon g_{ij}^1(x),}
\end{cases}\end{equation}
where
\begin{eqnarray}
\bar{H}_{ij}(\delta_{kl}+h_{kl},\frac{\partial h_{kl}}{\partial
t},\frac{\partial h_{kl}}{\partial
x^p})=H^{kl}\frac{\partial^2h_{ij}}{\partial x^kx^l}+d\frac{\partial
h_{ij}}{\partial
t}+\widetilde{H_{ij}}(\delta_{kl}+h_{kl},\frac{\partial
h_{kl}}{\partial t},\frac{\partial h_{kl}}{\partial x^p}).
\end{eqnarray}
By the definition (4.7) and (5.2)-(5.3), we have
\begin{eqnarray}
\nonumber &&\bar{H}_{ij}(\delta_{kl}+h_{kl},\frac{\partial
h_{kl}}{\partial t},\frac{\partial h_{kl}}{\partial
x^p})\vspace{2mm}\\ \nonumber&=&H^{kl}\frac{\partial^2h_{ij}}{\partial x^kx^l}\vspace{2mm}\\
\nonumber&&-\frac{1}{2}(\delta^{kl}+H^{kl})(\delta_{pq}+h_{pq})(\delta^{pa}+H^{pa})(\delta^{qb}+H^{qb})\vspace{2mm}\\
\nonumber&&\cdot\left(\frac{\partial h_{ai}}{\partial
x^k}+\frac{\partial h_{ak}}{\partial x^i}-\frac{\partial
h_{ik}}{\partial x^a}\right)\left(\frac{\partial h_{bj}}{\partial
x^l}+\frac{\partial h_{bl}}{\partial x^j}-\frac{\partial
h_{jl}}{\partial x^b}\right)\vspace{2mm}\\
\nonumber&&-\frac{1}{2}(\delta_{ik}+h_{ik})(\delta^{pr}+H^{pr})(\delta^{qs}+H^{qs})(\delta^{ka}+H^{ka})\vspace{2mm}\\
\nonumber&&\cdot\left(\frac{\partial h_{ar}}{\partial
x^s}+\frac{\partial h_{as}}{\partial x^r}-\frac{\partial
h_{rs}}{\partial x^a}\right)\left(\frac{\partial h_{pq}}{\partial
x^j}\right)\vspace{2mm}\\
\nonumber&&-\frac{1}{2}(\delta_{jk}+h_{jk})(\delta^{pr}+H^{pr})(\delta^{qs}+H^{qs})(\delta^{ka}+H^{ka})\vspace{2mm}\\
\nonumber&&\cdot\left(\frac{\partial h_{ar}}{\partial
x^s}+\frac{\partial h_{as}}{\partial x^r}-\frac{\partial
h_{rs}}{\partial x^a}\right)\left(\frac{\partial h_{pq}}{\partial
x^i}\right)\vspace{2mm}\\
\nonumber&&+2(\delta^{pq}+H^{pq})\frac{\partial h_{ip}}{\partial
t}\frac{\partial h_{jq}}{\partial
t}-2(\delta^{pq}+H^{pq})\frac{\partial h_{pq}}{\partial
t}\frac{\partial h_{ij}}{\partial t}\vspace{2mm}\\
\nonumber&&+\frac{1}{n-1}\left((\delta^{pq}+H^{pq})\frac{\partial
h_{pq}}{\partial
t}\right)^2(\delta_{ij}+h_{ij})-\frac{1}{n-1}(\delta^{pa}+H^{pa})(\delta^{qb}+H^{qb})\frac{\partial
h_{pq}}{\partial t}\frac{\partial h_{ab}}{\partial
t}(\delta_{ij}+h_{ij})\vspace{2mm}\\
\nonumber&=&-\frac{1}{2}\delta^{kl}\delta^{ab}\left(\frac{\partial
h_{ai}}{\partial x^k}+\frac{\partial h_{ak}}{\partial
x^i}-\frac{\partial h_{ik}}{\partial x^a}\right)\left(\frac{\partial
h_{bj}}{\partial x^l}+\frac{\partial h_{bl}}{\partial
x^j}-\frac{\partial h_{jl}}{\partial x^b}\right)\vspace{2mm}\\
\nonumber&&-\frac{1}{2}\delta^{pr}\delta^{qs}\left(\frac{\partial
h_{ir}}{\partial x^s}+\frac{\partial h_{is}}{\partial
x^r}-\frac{\partial h_{rs}}{\partial x^i}\right)\left(\frac{\partial
h_{pq}}{\partial x^j}\right) \vspace{2mm}\\
\nonumber&&-\frac{1}{2}\delta^{pr}\delta^{qs}\left(\frac{\partial
h_{jr}}{\partial x^s}+\frac{\partial h_{js}}{\partial
x^r}-\frac{\partial h_{rs}}{\partial x^j}\right)\left(\frac{\partial
h_{pq}}{\partial
x^i}\right)\vspace{2mm}\\
\nonumber&&+2\delta^{pq}\frac{\partial h_{ip}}{\partial
t}\frac{\partial h_{jq}}{\partial t}-2\delta^{pq}\frac{\partial
h_{pq}}{\partial
t}\frac{\partial h_{ij}}{\partial t}\vspace{2mm}\\
\nonumber&&+\frac{1}{n-1}\left(\delta^{pq}\frac{\partial
h_{pq}}{\partial
t}\right)^2\delta_{ij}-\frac{1}{n-1}\left(\delta^{pa}\delta^{qb}\frac{\partial
h_{pq}}{\partial t}\frac{\partial h_{ab}}{\partial
t}\right)\delta_{ij}-h^{kl}\frac{\partial^2h_{ij}}{\partial
x^kx^l}+O(||h_{kl}||+||Dh_{kl}||)^3\vspace{2mm}\\
&=&O(||h_{kl}||+||Dh_{kl}||+||\frac{\partial^2h_{ij}}{\partial
x^kx^l}||)^2+O(||h_{kl}||+||Dh_{kl}||)^3,
\end{eqnarray}
where $$Dh_{kl}\triangleq (\frac{\partial h_{kl}}{\partial
t},\frac{\partial h_{kl}}{\partial x^p})$$ and $||\cdot||$ stands
for the norm with respect to the flat metric $\delta_{ij}$.

By the theory of dissipative wave equations (see \cite{ma},
\cite{ni}), we know that, for sufficiently small $\epsilon>0$, the
Cauchy problem (5.5), i.e., (5.1), admits a unique smooth solution
for all $t\ge 0$ on $\mathbb{R}^n$ with $n\ge 3$. The proof of
Theorem 5.1 is completed. $\quad\Box$

\section{Wave character of curvatures --- Derivation of dissipative hyperbolic geometric flow}

In this section, we will illustrate why we choose (2.1) as the
equation of the dissipative hyperbolic geometric flow. Based on this
we derive the nonlinear wave equations satisfied by the curvatures.
The results presented in this section show the wave character of
curvatures.

We first assume that the metrics on a manifold $\mathscr{M}$ evolve
by the following equation
\begin{eqnarray}
\frac{\partial^2 g_{ij}}{\partial
t^2}(x,t)=-2R_{ij}(x,t)+G_{ij}(x,t),
\end{eqnarray}
where
\begin{eqnarray}
\nonumber{\displaystyle G_{ij}(x,t)}&=&{\displaystyle
ag^{pq}\frac{\partial g_{ip}}{\partial t}\frac{\partial
g_{jq}}{\partial t}+bg^{pq}\frac{\partial g_{pq}}{\partial
t}\frac{\partial g_{ij}}{\partial t}+d\frac{\partial
g_{ij}}{\partial t}+eg^{pq}\frac{\partial g_{pq}}{\partial
t}g_{ij}+}\vspace{2mm}\\& & {\displaystyle f(g^{pq}\frac{\partial
g_{pq}}{\partial t})^2g_{ij}+h\left(\frac{\partial g^{pq}}{\partial
t}\frac{\partial g_{pq}}{\partial t}\right)g_{ij}}
\end{eqnarray}
and the terms $a,\;b,\;d,\;e,\;f,\;h$ are all constants determined
below. Then direct calculation gives
\begin{eqnarray}
\nonumber{\displaystyle \Gamma^h_{jl}}&=& {\displaystyle
\frac{1}{2}g^{hm}(\frac{\partial g_{mj}}{\partial
x^l}+\frac{\partial g_{ml}}{\partial x^j}-\frac{\partial
g_{jl}}{\partial x^m}),}\vspace{2mm}\\ \nonumber{\displaystyle
\frac{\partial^2 R_{ik}}{\partial t^2}}&=&{\displaystyle
\frac{\partial^2 (g^{jl}R_{ijkl})}{\partial t^2}}\vspace{2mm}\\
\nonumber&=&{\displaystyle g^{jl}\frac{\partial^2 R_{ijkl}}{\partial
t^2}+2\frac{\partial g^{jl}}{\partial t}\frac{\partial
R_{ijkl}}{\partial t}+R_{ijkl}\frac{\partial^2 g^{jl}}{\partial
t^2},}\vspace{2mm}\\ \nonumber{\displaystyle \frac{\partial^2
R}{\partial t^2}}&=&{\displaystyle \frac{\partial^2
(g^{ik}R_{ik})}{\partial t^2}}\vspace{2mm}\\
\nonumber&=&{\displaystyle g^{ik}\frac{\partial^2 R_{ik}}{\partial
t^2}+2\frac{\partial g^{ik}}{\partial t}\frac{\partial
R_{ik}}{\partial t}+R_{ik}\frac{\partial^2 g^{ik}}{\partial t^2}}
\vspace{2mm}\\ \nonumber&=&{\displaystyle
g^{ik}\left(g^{jl}\frac{\partial^2 R_{ijkl}}{\partial
t^2}+2\frac{\partial g^{jl}}{\partial t}\frac{\partial
R_{ijkl}}{\partial t}+R_{ijkl}\frac{\partial^2 g^{jl}}{\partial
t^2}\right)+2\frac{\partial g^{ik}}{\partial t}\frac{\partial
R_{ik}}{\partial t}+R_{ik}\frac{\partial^2 g^{ik}}{\partial
t^2}}\vspace{2mm}\\ \nonumber&=&{\displaystyle
g^{ik}g^{jl}\frac{\partial^2 R_{ijkl}}{\partial t^2}+2\frac{\partial
g^{jl}}{\partial t}(\frac{\partial R_{jl}}{\partial
t}-\frac{\partial g^{ik}}{\partial t}R_{ijkl})+2\frac{\partial
g^{ik}}{\partial t}\frac{\partial R_{ik}}{\partial
t}+2R_{ik}\frac{\partial^2 g^{ik}}{\partial t^2}}\vspace{2mm}\\
\nonumber&=&{\displaystyle g^{ik}g^{jl}\frac{\partial^2
R_{ijkl}}{\partial t^2}+4\frac{\partial g^{ik}}{\partial
t}\frac{\partial R_{ik}}{\partial t}-2\frac{\partial
g^{jl}}{\partial t}\frac{\partial g^{ik}}{\partial
t}R_{ijkl}}\vspace{2mm}\\ \nonumber& & {\displaystyle
+2R_{ik}\left(-g^{ir}g^{ks}\frac{\partial^2 g_{rs}}{\partial
t^2}+2g^{ir}g^{ks}g^{pq}\frac{\partial g_{pr}}{\partial
t}\frac{\partial g_{qs}}{\partial t}\right)}\vspace{2mm}\\
\nonumber&=&{\displaystyle g^{ik}g^{jl}\frac{\partial^2
R_{ijkl}}{\partial t^2}+4\frac{\partial g^{ik}}{\partial
t}\frac{\partial R_{ik}}{\partial t}-2\frac{\partial
g^{jl}}{\partial t}\frac{\partial g^{ik}}{\partial
t}R_{ijkl}+}\vspace{2mm}\\ \nonumber& &{\displaystyle
4g^{ir}g^{ks}g^{pq}\frac{\partial g_{pr}}{\partial t}\frac{\partial
g_{qs}}{\partial
t}R_{ik}-2g^{ir}g^{ks}(-2R_{rs}+G_{rs})R_{ik}}\vspace{2mm}\\
\nonumber&= &{\displaystyle g^{ik}g^{jl}\frac{\partial^2
R_{ijkl}}{\partial t^2}+4\frac{\partial g^{ik}}{\partial
t}\frac{\partial R_{ik}}{\partial t}-2\frac{\partial
g^{jl}}{\partial t}\frac{\partial g^{ik}}{\partial
t}R_{ijkl}+}\vspace{2mm}\\& & {\displaystyle
4g^{ir}g^{ks}g^{pq}\frac{\partial g_{pr}}{\partial t}\frac{\partial
g_{qs}}{\partial t}R_{ik}+4|Ric|^2-2g^{ir}g^{ks}R_{ik}G_{rs},}
\end{eqnarray}
where
$$|Ric|^2=g^{ik}g^{jl}R_{ij}R_{kl}$$ is the norm of Ricci curvature
tensor $Ric=R_{ik}$. In (6.3), we have made use of the evolution
equation (6.1).

We choose the normal coordinates around a fixed point $p$ on the
manifold $\mathscr{M}$ such that
$$\Gamma_{ij}^k(p)=0\;\;(\;\;\forall\;\; i,\;j,\;k), \quad {\rm or
\,
equivalently,}\quad \frac{\partial g_{ij}}{\partial x^k}(p)=0,$$
where $\Gamma_{ij}^k$ stand for the Christoffel symbols. By the
computations (5.2)-(5.5) in \cite{dkl}, we get
\begin{eqnarray}
\nonumber\dfrac{\partial^{2}}{\partial
t^{2}}R_{ijkl}&=&\dfrac{1}{2}\left[\dfrac{\partial^{2}}{\partial
x^{i}\partial x^{l}}\left(\dfrac{\partial^{2} g_{kj}}{\partial
t^{2}}\right)+\dfrac{\partial^{2}}{\partial x^{i}\partial
x^{j}}\left(\dfrac{\partial^{2} g_{kl}}{\partial t^{2}}\right)-
\dfrac{\partial^{2}}{\partial x^{i}\partial
x^{k}}\left(\dfrac{\partial^{2} g_{jl}}{\partial
t^{2}}\right)\right]\vspace{2mm}\\
\nonumber&&-\dfrac{1}{2}\left[\dfrac{\partial^{2}}{\partial
x^{j}\partial x^{l}}\left(\dfrac{\partial^{2} g_{ki}}{\partial
t^{2}}\right)+\dfrac{\partial^{2}}{\partial x^{j}\partial
x^{i}}\left(\dfrac{\partial^{2} g_{kl}}{\partial
t^{2}}\right)-\dfrac{\partial^{2}}{\partial x^{j}\partial
x^{k}}\left(\dfrac{\partial^{2} g_{il}}{\partial
t^{2}}\right)\right]\vspace{2mm}\\&&+2g_{pq}\left(\dfrac{\partial}{\partial
t}\Gamma^{p}_{il}\cdot\dfrac{\partial}{\partial
t}\Gamma^{q}_{jk}-\dfrac{\partial}{\partial
t}\Gamma^{p}_{jl}\cdot\dfrac{\partial}{\partial
t}\Gamma^{q}_{ik}\right).
\end{eqnarray}
Then it follows from (6.1) and (6.4) that
\begin{eqnarray}
\nonumber\dfrac{\partial^{2}}{\partial
t^{2}}R_{ijkl}&=&\dfrac{1}{2}\left[\dfrac{\partial^{2}}{\partial
x^{i}\partial x^{l}}(-2R_{kj})+\dfrac{\partial^{2}}{\partial
x^{i}\partial x^{j}}(-2R_{kl})- \dfrac{\partial^{2}}{\partial
x^{i}\partial x^{k}}(-2R_{jl})\right]\vspace{2mm}\\
\nonumber &&-\dfrac{1}{2}\left[\dfrac{\partial^{2}}{\partial
x^{j}\partial x^{l}}(-2R_{ki})+\dfrac{\partial^{2}}{\partial
x^{j}\partial x^{i}}(-2R_{kl})-\dfrac{\partial^{2}}{\partial
x^{j}\partial x^{k}}(-2R_{il})\right]\vspace{2mm}\\
\nonumber&&+2g_{pq}\left(\dfrac{\partial}{\partial
t}\Gamma^{p}_{il}\cdot\dfrac{\partial}{\partial
t}\Gamma^{q}_{jk}-\dfrac{\partial}{\partial
t}\Gamma^{p}_{jl}\cdot\dfrac{\partial}{\partial
t}\Gamma^{q}_{ik}\right)\vspace{2mm}\\
\nonumber&&+\dfrac{1}{2}\left[\dfrac{\partial^{2}}{\partial
x^{i}\partial x^{l}}G_{kj}+\dfrac{\partial^{2}}{\partial
x^{i}\partial x^{j}}G_{kl}- \dfrac{\partial^{2}}{\partial
x^{i}\partial
x^{k}}G_{jl}\right]\vspace{2mm}\\&&-\dfrac{1}{2}\left[\dfrac{\partial^{2}}{\partial
x^{j}\partial x^{l}}G_{ki}+\dfrac{\partial^{2}}{\partial
x^{j}\partial x^{i}}G_{kl}-\dfrac{\partial^{2}}{\partial
x^{j}\partial x^{k}}G_{il}\right].
\end{eqnarray}
Similar to Hamilton \cite{ha}, by Theorem 5.1 in \cite{dkl} we have
\begin{eqnarray}
\nonumber&&\;\;\dfrac{1}{2}\left[\dfrac{\partial^{2}}{\partial
x^{i}\partial x^{l}}(-2R_{kj})+\dfrac{\partial^{2}}{\partial
x^{i}\partial x^{j}}(-2R_{kl})- \dfrac{\partial^{2}}{\partial
x^{i}\partial x^{k}}(-2R_{jl})\right]\vspace{2mm}\\
\nonumber&&-\dfrac{1}{2}\left[\dfrac{\partial^{2}}{\partial
x^{j}\partial x^{l}}(-2R_{ki})+\dfrac{\partial^{2}}{\partial
x^{j}\partial x^{i}}(-2R_{kl})-\dfrac{\partial^{2}}{\partial
x^{j}\partial x^{k}}(-2R_{il})\right]\vspace{2mm}\\
\nonumber&&+2g_{pq}\left(\dfrac{\partial}{\partial
t}\Gamma^{p}_{il}\cdot\dfrac{\partial}{\partial
t}\Gamma^{q}_{jk}-\dfrac{\partial}{\partial
t}\Gamma^{p}_{jl}\cdot\dfrac{\partial}{\partial
t}\Gamma^{q}_{ik}\right)\vspace{2mm}\\ \nonumber&=& \triangle
R_{ijkl}+2\left(B_{ijkl}-B_{ijlk}-B_{iljk}+B_{ikjl}\right)\vspace{2mm}\\
\nonumber & &
-g^{pq}\left(R_{pjkl}R_{qi}+R_{ipkl}R_{qj}+R_{ijpl}R_{qk}+R_{ijkp}R_{ql}\right)\vspace{2mm}\\
& &+2g_{pq}\left(\dfrac{\partial}{\partial
t}\Gamma^{p}_{il}\cdot\dfrac{\partial}{\partial
t}\Gamma^{q}_{jk}-\dfrac{\partial}{\partial
t}\Gamma^{p}_{jl}\cdot\dfrac{\partial}{\partial
t}\Gamma^{q}_{ik}\right),
\end{eqnarray}
where $B_{ijkl}=g^{pr}g^{qs}R_{piqj}R_{rksl}$ and $\triangle$ is the
Laplacian with respect to the evolving metric.

Combing (6.3), (6.5) and (6.6) and referring to the computations in
Theorem 5.3 in \cite{dkl} leads to
\begin{eqnarray}
\nonumber\frac{\partial^2 R}{\partial t^2}&=&\Delta R+2|{\rm
Ric}|^2\vspace{2mm}\\ \nonumber&&
+2g^{ik}g^{jl}g_{pq}\left(\dfrac{\partial}{\partial
t}\Gamma^{p}_{il}\dfrac{\partial}{\partial
t}\Gamma^{q}_{jk}-\dfrac{\partial}{\partial
t}\Gamma^{p}_{jl}\dfrac{\partial}{\partial
t}\Gamma^{q}_{ik}\right)\vspace{2mm}\\ \nonumber&
&-2g^{ik}g^{jp}g^{lq}\dfrac{\partial g_{pq}}{\partial
t}\dfrac{\partial}{\partial t}R_{ijkl}\vspace{2mm}\\ \nonumber&
&-2g^{ip}g^{kq}\dfrac{\partial g_{pq}}{\partial t}\dfrac{\partial
R_{ik}}{\partial t}+4R_{ik}g^{ip}g^{rq}g^{sk}\dfrac{\partial
g_{pq}}{\partial t}\dfrac{\partial g_{rs}}{\partial t}\vspace{2mm}\\
\nonumber
&&+\dfrac{1}{2}g^{ik}g^{jl}\left[\dfrac{\partial^{2}}{\partial
x^{i}\partial x^{l}}G_{kj}+\dfrac{\partial^{2}}{\partial
x^{i}\partial x^{j}}G_{kl}- \dfrac{\partial^{2}}{\partial
x^{i}\partial x^{k}}G_{jl}\right]\vspace{2mm}\\
\nonumber&&-\dfrac{1}{2}g^{ik}g^{jl}\left[\dfrac{\partial^{2}}{\partial
x^{j}\partial x^{l}}G_{ki}+\dfrac{\partial^{2}}{\partial
x^{j}\partial x^{i}}G_{kl}-\dfrac{\partial^{2}}{\partial
x^{j}\partial x^{k}}G_{il}\right]\vspace{2mm}\\
\nonumber&&-2g^{ir}g^{ks}R_{ik}G_{rs} \vspace{2mm}\\
\nonumber&=&\Delta R+2|{\rm Ric}|^2\vspace{2mm}\\ \nonumber&&
+2g^{ik}g^{jl}g_{pq}\left(\dfrac{\partial}{\partial
t}\Gamma^{p}_{il}\dfrac{\partial}{\partial
t}\Gamma^{q}_{jk}-\dfrac{\partial}{\partial
t}\Gamma^{p}_{jl}\dfrac{\partial}{\partial
t}\Gamma^{q}_{ik}\right)\vspace{2mm}\\ \nonumber&
&-2g^{ik}g^{jp}g^{lq}\dfrac{\partial g_{pq}}{\partial
t}\dfrac{\partial}{\partial t}R_{ijkl}\vspace{2mm}\\ \nonumber&
&-2g^{ip}g^{kq}\dfrac{\partial g_{pq}}{\partial t}\dfrac{\partial
R_{ik}}{\partial t}+4R_{ik}g^{ip}g^{rq}g^{sk}\dfrac{\partial
g_{pq}}{\partial t}\dfrac{\partial g_{rs}}{\partial t}\vspace{2mm}\\
&&+g^{ik}g^{jl}\left(\dfrac{\partial^{2}}{\partial x^{i}\partial
x^{l}}G_{kj}-\dfrac{\partial^{2}}{\partial x^{i}\partial
x^{k}}G_{jl}\right)-2g^{ir}g^{ks}R_{ik}G_{rs}.
\end{eqnarray}

In the normal coordinates, we have
\begin{eqnarray}
\nonumber & &g^{ik}g^{jl}\left(\dfrac{\partial^{2}}{\partial
x^{i}\partial x^{l}}G_{kj}-\dfrac{\partial^{2}}{\partial
x^{i}\partial x^{k}}G_{jl}\right)-2g^{ir}g^{ks}R_{ik}G_{rs}\vspace{2mm}\\
\nonumber&=&g^{ik}g^{jl}\left(\nabla_i\nabla_l G_{kj}
+\nabla_i\Gamma^p_{lk}G_{pj}+\nabla_i\Gamma^p_{lj}G_{pk}\right)\vspace{2mm}\\
\nonumber&&-g^{ik}g^{jl}\left(\nabla_j\nabla_l G_{ki}
+\nabla_j\Gamma^p_{lk}G_{pi}+\nabla_j\Gamma^p_{li}G_{pk}\right)\vspace{2mm}\\
\nonumber&&-2g^{ir}g^{ks}R_{ik}G_{rs}\vspace{2mm}\\ \nonumber&=&
g^{ik}g^{jl}(\nabla_i\nabla_l G_{kj}-\nabla_j\nabla_l
G_{ki})+g^{ik}g^{jl}\nabla_i\Gamma^p_{lk}G_{pj}-g^{ik}g^{jl}\nabla_j\Gamma^p_{lk}G_{pi}
\vspace{2mm}\\
\nonumber&&+g^{ik}g^{jl}(\nabla_i\Gamma^p_{lj}-\nabla_j\Gamma^p_{li})G_{pk}-2g^{ir}g^{ks}R_{ik}G_{rs}
\vspace{2mm}\\ \nonumber&=&g^{ik}g^{jl}(\nabla_i\nabla_l
G_{kj}-\nabla_j\nabla_l
G_{ki})+g^{ik}g^{jl}R_{ijl}^pG_{pk}-2g^{ir}g^{ks}R_{ik}G_{rs}
\vspace{2mm}\\&=&g^{ik}g^{jl}(\nabla_i\nabla_l
G_{kj}-\nabla_j\nabla_l G_{ki})-g^{ir}g^{ks}R_{ik}G_{rs},
\end{eqnarray}
where we have made use of the following equality in the normal
coordinates
\begin{eqnarray*}
R_{ijl}^p=\frac{\partial \Gamma^p_{lj}}{\partial x^i}-\frac{\partial
\Gamma^p_{li}}{\partial x^j}
=\nabla_i\Gamma^p_{lj}-\nabla_j\Gamma^p_{li}
\end{eqnarray*}
and $\nabla_l$ means the covariant derivative in the direction
$\frac{\partial }{\partial x^l}$. In the normal coordinates, we
easily obtain
\begin{eqnarray*}
\nonumber\frac{\partial g_{ij}}{\partial
x^k}&=&g_{ip}\Gamma_{kj}^p+g_{jp}\Gamma_{ik}^p,\vspace{2mm}\\
\nonumber \frac{\partial ^2 g_{ij}}{\partial x^k\partial
t}&=&g_{ip}\frac{\partial \Gamma_{kj}^p}{\partial
t}+g_{jp}\frac{\partial \Gamma_{ki}^p}{\partial t},\vspace{2mm}\\
\nonumber \frac{\partial ^2 g_{ij}}{\partial x^k\partial
x^l}&=&g_{ip}\frac{\partial \Gamma_{lj}^p}{\partial
x^k}+g_{jp}\frac{\partial \Gamma_{li}^p}{\partial x^k}.
\end{eqnarray*}
Then, we have
\begin{eqnarray*}
\dfrac{\partial }{\partial t} \left(\frac{\partial ^2
g_{ij}}{\partial x^k\partial x^l}\right)=\frac{\partial
g_{ip}}{\partial t}\frac{\partial \Gamma_{lj}^p}{\partial
x^k}+\frac{\partial g_{jp}}{\partial t}\frac{\partial
\Gamma_{li}^p}{\partial x^k}+g_{ip}\dfrac{\partial }{\partial
t}\frac{\partial \Gamma_{lj}^p}{\partial x^k}+g_{jp}\dfrac{\partial
}{\partial t}\frac{\partial \Gamma_{li}^p}{\partial x^k}.
\end{eqnarray*}
This implies that
\begin{eqnarray*}
\nonumber\nabla_l(\frac{\partial g_{jp}}{\partial
t})&=&\frac{\partial^2 g_{jp}}{\partial x^l\partial
t}-\Gamma_{lj}^q\frac{\partial g_{pq}}{\partial
t}-\Gamma_{lp}^q\frac{\partial g_{jq}}{\partial t}=\frac{\partial^2
g_{ip}}{\partial x^l\partial t},\vspace{2mm}\\
\nonumber\nabla_i\nabla_l(\frac{\partial g_{jp}}{\partial
t})&=&\dfrac{\partial }{\partial x^i}(\frac{\partial^2
g_{jp}}{\partial x^l\partial t}-\Gamma_{lj}^q\frac{\partial
g_{pq}}{\partial t}-\Gamma_{lp}^q\frac{\partial g_{jq}}{\partial
t})-\Gamma_{il}^r(\frac{\partial^2 g_{jp}}{\partial x^r\partial
t}-\Gamma_{rj}^q\frac{\partial g_{pq}}{\partial
t}-\Gamma_{rp}^q\frac{\partial g_{jq}}{\partial t})\vspace{2mm}\\
\nonumber& &-\Gamma_{ij}^r(\frac{\partial^2 g_{rp}}{\partial
x^l\partial t}-\Gamma_{lr}^q\frac{\partial g_{pq}}{\partial
t}-\Gamma_{lp}^q\frac{\partial g_{rq}}{\partial
t})-\Gamma_{ip}^r(\frac{\partial^2 g_{jr}}{\partial x^l\partial
t}-\Gamma_{lj}^q\frac{\partial g_{rq}}{\partial
t}-\Gamma_{rl}^q\frac{\partial g_{jq}}{\partial t})\vspace{2mm}\\
\nonumber &=&\dfrac{\partial }{\partial t} (\frac{\partial ^2
g_{ij}}{\partial x^k\partial x^l})-\dfrac{\partial
\Gamma_{lj}^q}{\partial x^i}\frac{\partial g_{pq}}{\partial
t}-\dfrac{\partial \Gamma_{lp}^q}{\partial x^i}\frac{\partial
g_{jq}}{\partial t}\vspace{2mm}\\ \nonumber&=&g_{jr}\dfrac{\partial
}{\partial t}\frac{\partial \Gamma_{lp}^r}{\partial
x^i}+g_{pr}\dfrac{\partial }{\partial t}\frac{\partial
\Gamma_{lj}^r}{\partial x^i}.
\end{eqnarray*}
By the   direct computations, we have
\begin{eqnarray}
\nonumber& &g^{ik}g^{jl}\nabla_i\nabla_l(g^{pq}\frac{\partial
g_{jp}}{\partial t}\frac{\partial g_{kq}}{\partial
t})-g^{ik}g^{jl}\nabla_j\nabla_l(g^{pq}\frac{\partial
g_{ip}}{\partial t}\frac{\partial g_{kq}}{\partial t})\vspace{2mm}\\
\nonumber &=&g^{ik}g^{jl}g^{pq}\nabla_i\nabla_l(\frac{\partial
g_{jp}}{\partial t}\frac{\partial g_{kq}}{\partial
t})-g^{ik}g^{jl}g^{pq}\nabla_j\nabla_l(\frac{\partial
g_{ip}}{\partial t}\frac{\partial g_{kq}}{\partial t})\vspace{2mm}\\
\nonumber &=&g^{ik}g^{jl}g^{pq}\left[\nabla_i\nabla_l(\frac{\partial
g_{jp}}{\partial t})\frac{\partial g_{kq}}{\partial
t}+\frac{\partial g_{jp}}{\partial t}\nabla_i\nabla_l(\frac{\partial
g_{kq}}{\partial t})\right]\vspace{2mm}\\ \nonumber& &
-g^{ik}g^{jl}g^{pq}\left[\nabla_j\nabla_l(\frac{\partial
g_{ip}}{\partial t})\frac{\partial g_{kq}}{\partial
t}+\frac{\partial g_{ip}}{\partial t}\nabla_j\nabla_l(\frac{\partial
g_{kq}}{\partial t})\right]\vspace{2mm}\\ \nonumber& &
+g^{ik}g^{jl}g^{pq}\left[\nabla_i(\frac{\partial g_{jp}}{\partial
t})\nabla_l(\frac{\partial g_{kq}}{\partial
t})+\nabla_l(\frac{\partial g_{jp}}{\partial
t})\nabla_i(\frac{\partial g_{kq}}{\partial t})\right]\vspace{2mm}\\
\nonumber& &-g^{ik}g^{jl}g^{pq}\left[\nabla_j(\frac{\partial
g_{ip}}{\partial t})\nabla_l(\frac{\partial g_{kq}}{\partial
t})+\nabla_l(\frac{\partial g_{ip}}{\partial
t})\nabla_j(\frac{\partial g_{kq}}{\partial t})\right]\vspace{2mm}\\
\nonumber&= &g^{ik}g^{jl}g^{pq}\left[\frac{\partial g_{jp}}{\partial
t}g_{kr}\dfrac{\partial }{\partial t}\frac{\partial
\Gamma_{lq}^r}{\partial x^i}+\frac{\partial g_{jp}}{\partial
t}g_{qr}\dfrac{\partial }{\partial t}\frac{\partial
\Gamma_{kl}^r}{\partial x^i}+\frac{\partial g_{kq}}{\partial
t}g_{jr}\dfrac{\partial }{\partial t}\frac{\partial
\Gamma_{lp}^r}{\partial x^i}+\frac{\partial g_{kq}}{\partial
t}g_{pr}\dfrac{\partial }{\partial t}\frac{\partial
\Gamma_{lj}^r}{\partial x^i}\right]\vspace{2mm}\\ \nonumber&
&-g^{ik}g^{jl}g^{pq}\left[\frac{\partial g_{ip}}{\partial
t}g_{kr}\dfrac{\partial }{\partial t}\frac{\partial
\Gamma_{lq}^r}{\partial x^j}+\frac{\partial g_{ip}}{\partial
t}g_{qr}\dfrac{\partial }{\partial t}\frac{\partial
\Gamma_{kl}^r}{\partial x^j}+\frac{\partial g_{kq}}{\partial
t}g_{ir}\dfrac{\partial }{\partial t}\frac{\partial
\Gamma_{lp}^r}{\partial x^j}+\frac{\partial g_{kq}}{\partial
t}g_{pr}\dfrac{\partial }{\partial t}\frac{\partial
\Gamma_{li}^r}{\partial x^j}\right]\vspace{2mm}\\ \nonumber&
&+g^{ik}g^{jl}g^{pq}\left[\frac{\partial^2 g_{jp}}{\partial
t\partial x^i}\frac{\partial^2 g_{kq}}{\partial t\partial
x^l}+\frac{\partial^2 g_{jp}}{\partial t\partial
x^l}\frac{\partial^2 g_{kq}}{\partial t\partial
x^i}-\frac{\partial^2 g_{ip}}{\partial t\partial
x^j}\frac{\partial^2 g_{kq}}{\partial t\partial
x^l}-\frac{\partial^2 g_{ip}}{\partial t\partial
x^l}\frac{\partial^2 g_{kq}}{\partial t\partial
x^j}\right]\vspace{2mm}\\
\nonumber&=&\left(g^{jl}g^{pq}\frac{\partial g_{jp}}{\partial
t}\dfrac{\partial }{\partial t}\frac{\partial
\Gamma_{lq}^i}{\partial x^i}+g^{ik}g^{jl}\frac{\partial
g_{jp}}{\partial t}\dfrac{\partial }{\partial t}\frac{\partial
\Gamma_{kl}^p}{\partial x^i}+g^{ik}g^{pq}\frac{\partial
g_{kq}}{\partial t}\dfrac{\partial }{\partial t}\frac{\partial
\Gamma_{lp}^l}{\partial x^i}+g^{ik}g^{jl}\frac{\partial
g_{kq}}{\partial t}\dfrac{\partial }{\partial t}\frac{\partial
\Gamma_{lj}^q}{\partial x^i}\right)\vspace{2mm}\\
\nonumber&&-\left(g^{jl}g^{pq}\frac{\partial g_{ip}}{\partial
t}\dfrac{\partial }{\partial t}\frac{\partial
\Gamma_{lq}^i}{\partial x^j}+g^{ik}g^{jl}\frac{\partial
g_{ip}}{\partial t}\dfrac{\partial }{\partial t}\frac{\partial
\Gamma_{kl}^p}{\partial x^j}+g^{jl}g^{pq}\frac{\partial
g_{kq}}{\partial t}\dfrac{\partial }{\partial t}\frac{\partial
\Gamma_{lp}^k}{\partial x^j}+g^{ik}g^{jl}\frac{\partial
g_{kq}}{\partial t}\dfrac{\partial }{\partial t}\frac{\partial
\Gamma_{li}^q}{\partial x^j}\right)\vspace{2mm}\\
\nonumber&&+g^{ik}g^{jl}g^{pq}\left[\frac{\partial^2
g_{jp}}{\partial t\partial x^i}\frac{\partial^2 g_{kq}}{\partial
t\partial x^l}+\frac{\partial^2 g_{jp}}{\partial t\partial
x^l}\frac{\partial^2 g_{kq}}{\partial t\partial
x^i}-\frac{\partial^2 g_{ip}}{\partial t\partial
x^j}\frac{\partial^2 g_{kq}}{\partial t\partial
x^l}-\frac{\partial^2 g_{ip}}{\partial t\partial
x^l}\frac{\partial^2 g_{kq}}{\partial t\partial
x^j}\right]\vspace{2mm}\\
\nonumber&=&\left(g^{ik}g^{jl}\frac{\partial g_{kl}}{\partial
t}\frac{\partial }{\partial t}R_{ij}+2g^{ik}g^{jl}\frac{\partial
g_{kl}}{\partial t}\frac{\partial }{\partial t}\frac{\partial
\Gamma_{jp}^p}{\partial x^i}+g^{ik}g^{jl}\frac{\partial
g_{kp}}{\partial t}\frac{\partial }{\partial
t}R_{ijl}^p-2g^{ik}g^{jl}\frac{\partial g_{ip}}{\partial
t}\frac{\partial }{\partial t}\frac{\partial \Gamma_{kl}^p}{\partial
x^j}\right)\vspace{2mm}\\ \nonumber&
&+g^{ik}g^{jl}g^{pq}\left(\frac{\partial^2 g_{jp}}{\partial
t\partial x^i}\frac{\partial^2 g_{kq}}{\partial t\partial
x^l}+\frac{\partial^2 g_{jp}}{\partial t\partial
x^l}\frac{\partial^2 g_{kq}}{\partial t\partial
x^i}-2\frac{\partial^2 g_{ip}}{\partial t\partial
x^j}\frac{\partial^2 g_{kq}}{\partial t\partial x^l}\right)\vspace{2mm}\\
\nonumber&=&\left(-2\frac{\partial g^{ij} }{\partial
t}\frac{\partial R_{ij}}{\partial t}-2\frac{\partial g^{ij}
}{\partial t}\frac{\partial }{\partial t}\dfrac{\partial
\Gamma_{jp}^p}{\partial x^i}-g_{kl}\frac{\partial g^{ik}}{\partial
t}\frac{\partial g^{jl}}{\partial t}R_{ij}+\frac{\partial
g^{ik}}{\partial t}\frac{\partial g^{jl}}{\partial
t}R_{ijkl}-2g^{ik}g^{jl}\frac{\partial g_{ip}}{\partial
t}\frac{\partial }{\partial t}\frac{\partial \Gamma_{kl}^p}{\partial
x^j}\right)\vspace{2mm}\\ &
&+g^{ik}g^{jl}g^{pq}\left(\frac{\partial^2 g_{jp}}{\partial
t\partial x^i}\frac{\partial^2 g_{kq}}{\partial t\partial
x^l}+\frac{\partial^2 g_{jp}}{\partial t\partial
x^l}\frac{\partial^2 g_{kq}}{\partial t\partial
x^i}-2\frac{\partial^2 g_{ip}}{\partial t\partial
x^j}\frac{\partial^2 g_{kq}}{\partial t\partial x^l}\right),
\end{eqnarray}
where we have made use of the following equation
\begin{eqnarray*}
g^{jl}g^{pq}\frac{\partial g_{jp}}{\partial t}\dfrac{\partial
}{\partial t}\frac{\partial \Gamma_{lq}^i}{\partial
x^i}&=&g^{ik}g^{jl}\frac{\partial g_{kl}}{\partial t}\dfrac{\partial
}{\partial t}\frac{\partial \Gamma_{ij}^p}{\partial
x^p}\vspace{2mm}\\&=&g^{ik}g^{jl}\frac{\partial g_{kl}}{\partial
t}\dfrac{\partial }{\partial t}(R_{pij}^p+\frac{\partial
\Gamma_{jp}^p}{\partial
x^i})\vspace{2mm}\\&=&g^{ik}g^{jl}\frac{\partial g_{kl}}{\partial
t}\dfrac{\partial }{\partial t}R_{ij}+g^{ik}g^{jl}\frac{\partial
g_{kl}}{\partial t}\dfrac{\partial }{\partial t}\frac{\partial
\Gamma_{jp}^p}{\partial x^i}.
\end{eqnarray*}

Analogously, we obtain
\begin{eqnarray} \nonumber\quad&
&g^{ik}g^{jl}\nabla_i\nabla_l(g^{pq}\frac{\partial g_{pq}}{\partial
t}\frac{\partial g_{jk}}{\partial
t})-g^{ik}g^{jl}\nabla_j\nabla_l(g^{pq}\frac{\partial
g_{pq}}{\partial t}\frac{\partial g_{ik}}{\partial t})\vspace{2mm}\\
\nonumber &=&g^{ik}g^{jl}g^{pq}\nabla_i\nabla_l(\frac{\partial
g_{pq}}{\partial t}\frac{\partial g_{jk}}{\partial
t})-g^{ik}g^{jl}g^{pq}\nabla_j\nabla_l(\frac{\partial
g_{pq}}{\partial t}\frac{\partial g_{ik}}{\partial t})\vspace{2mm}\\
\nonumber &=&(g^{pq}\frac{\partial g_{pq}}{\partial
t})(\frac{\partial R}{\partial t}-\frac{\partial g^{jl}}{\partial
t}R_{jl})-2(g^{pq}\frac{\partial g_{pq}}{\partial
t})g^{ik}\dfrac{\partial }{\partial t}\frac{\partial
\Gamma_{kr}^r}{\partial x^i}-2\frac{\partial g^{ik}}{\partial
t}\dfrac{\partial }{\partial t}\frac{\partial
\Gamma_{kr}^r}{\partial x^i}\vspace{2mm}\\&
&+g^{ik}g^{jl}g^{pq}\left(\dfrac{\partial^2 g_{pq}}{\partial
t\partial x^i}\dfrac{\partial^2 g_{jk}}{\partial t\partial
x^l}+\dfrac{\partial^2 g_{pq}}{\partial t\partial
x^l}\dfrac{\partial^2 g_{jk}}{\partial t\partial
x^i}-2\dfrac{\partial^2 g_{pq}}{\partial t\partial
x^j}\dfrac{\partial^2 g_{ik}}{\partial t\partial x^l}\right).
\end{eqnarray}

On the other hand,
\begin{eqnarray}\nonumber&
&g^{ik}g^{jl}\nabla_i\nabla_l(\frac{\partial g_{jk}}{\partial
t})-g^{ik}g^{jl}\nabla_j\nabla_l(\frac{\partial g_{ik}}{\partial
t})\vspace{2mm}\\
\nonumber&=&g^{ik}g^{jl}\left(g_{jp}\dfrac{\partial }{\partial
t}\frac{\partial \Gamma_{kl}^p}{\partial x^i}+g_{kp}\dfrac{\partial
}{\partial t}\frac{\partial \Gamma_{jl}^p}{\partial
x^i}\right)-g^{ik}g^{jl}\left(g_{ip}\dfrac{\partial }{\partial
t}\frac{\partial \Gamma_{kl}^p}{\partial x^j}+g_{kp}\dfrac{\partial
}{\partial t}\frac{\partial \Gamma_{il}^p}{\partial
x^j}\right)\vspace{2mm}\\ \nonumber&=&g^{ik}\dfrac{\partial
}{\partial t}\frac{\partial \Gamma_{kl}^l}{\partial
x^i}+g^{jl}\dfrac{\partial }{\partial t}\frac{\partial
\Gamma_{jl}^i}{\partial x^i}-g^{jl}\dfrac{\partial }{\partial
t}\frac{\partial \Gamma_{lk}^k}{\partial x^j}-g^{jl}\dfrac{\partial
}{\partial t}\frac{\partial \Gamma_{li}^i}{\partial
x^j}\vspace{2mm}\\ \nonumber&=&g^{jl}\dfrac{\partial }{\partial
t}\frac{\partial \Gamma_{jl}^r}{\partial x^r}-g^{jl}\dfrac{\partial
}{\partial t}\frac{\partial \Gamma_{lr}^r}{\partial
x^j}\vspace{2mm}\\ &=&g^{jl}\dfrac{\partial }{\partial t}R_{jl},
\vspace{5mm}\\ \nonumber &
&g^{ik}g^{jl}\nabla_i\nabla_l\left(g^{pq}\frac{\partial
g_{pq}}{\partial
t}g_{jk}\right)-g^{ik}g^{jl}\nabla_j\nabla_l\left(g^{pq}\frac{\partial
g_{pq}}{\partial t}g_{ik}\right)\vspace{2mm}\\ \nonumber
&=&g^{ik}g^{jl}g^{pq}g_{kj}\nabla_i\nabla_l(\frac{\partial
g_{pq}}{\partial
t})-g^{ik}g^{jl}g^{pq}g_{ik}\nabla_j\nabla_l(\frac{\partial
g_{pq}}{\partial t})\vspace{2mm}\\ \nonumber
&=&(1-n)g^{ik}g^{jl}\nabla_i\nabla_k(\frac{\partial g_{jl}}{\partial
t})\vspace{2mm}\\
\nonumber&=&(1-n)g^{ik}g^{jl}\left(g_{jp}\dfrac{\partial }{\partial
t}\frac{\partial \Gamma_{kl}^p}{\partial x^i}+g_{lp}\dfrac{\partial
}{\partial t}\frac{\partial \Gamma_{kj}^p}{\partial
x^i}\right)\vspace{2mm}\\ &=&-2(n-1)g^{ik}\dfrac{\partial }{\partial
t}\frac{\partial \Gamma_{kl}^l}{\partial x^i}, \vspace{5mm}\\
\nonumber & &g^{ik}g^{jl}\nabla_i\nabla_l\left((g^{pq}\frac{\partial
g_{pq}}{\partial
t})^2g_{jk}\right)-g^{ik}g^{jl}\nabla_j\nabla_l\left((g^{pq}\frac{\partial
g_{pq}}{\partial t})^2g_{ik}\right)\vspace{2mm}\\ \nonumber
&=&g^{ik}g^{jl}g_{kj}g^{pq}g^{rs}\nabla_i\nabla_l\left(\frac{\partial
g_{pq}}{\partial t}\frac{\partial g_{rs}}{\partial
t}\right)-g^{ik}g^{jl}g_{ki}g^{pq}g^{rs}\nabla_j\nabla_l\left(\frac{\partial
g_{pq}}{\partial t}\frac{\partial g_{rs}}{\partial
t}\right)\vspace{2mm}\\
\nonumber&=&(1-n)g^{jl}g^{pq}g^{rs}\nabla_j\nabla_l\left(\frac{\partial
g_{pq}}{\partial t}\frac{\partial g_{rs}}{\partial
t}\right)\vspace{2mm}\\
\nonumber&=&-(n-1)g^{jl}g^{pq}g^{rs}\left[\nabla_j\nabla_l(\frac{\partial
g_{pq}}{\partial t})\frac{\partial g_{rs}}{\partial
t}+\nabla_j\nabla_l(\frac{\partial g_{rs}}{\partial
t})\frac{\partial g_{pq}}{\partial t}+2\nabla_j(\frac{\partial
g_{pq}}{\partial t})\nabla_l(\frac{\partial g_{rs}}{\partial
t})\right]\vspace{2mm}\\ \nonumber
&=&-2(n-1)g^{jl}g^{pq}g^{rs}\left[\frac{\partial g_{rs}}{\partial
t}\left(g_{pi}\dfrac{\partial }{\partial t}\frac{\partial
\Gamma_{lq}^i}{\partial x^j}+g_{qi}\dfrac{\partial }{\partial
t}\frac{\partial \Gamma_{lp}^i}{\partial
x^j}\right)+\frac{\partial^2 g_{pq}}{\partial t\partial
x^j}\frac{\partial^2 g_{rs}}{\partial
t\partial x^l}\right]\vspace{2mm}\\
&=&-(n-1)\left[4(g^{pq}\frac{\partial g_{pq}}{\partial
t})g^{jl}\dfrac{\partial }{\partial t}\frac{\partial
\Gamma_{ls}^s}{\partial x^j}+2g^{jl}g^{pq}g^{rs}\frac{\partial^2
g_{pq}}{\partial t\partial x^j}\frac{\partial^2 g_{rs}}{\partial
t\partial x^l}\right]
\end{eqnarray}
and \begin{eqnarray}
\nonumber&&g^{ik}g^{jl}\nabla_i\nabla_l(\frac{\partial
g^{pq}}{\partial t}\frac{\partial g_{pq}}{\partial
t}g_{jk})-g^{ik}g^{jl}\nabla_j\nabla_l(\frac{\partial
g^{pq}}{\partial t}\frac{\partial g_{pq}}{\partial
t}g_{ik})\vspace{2mm}\\ \nonumber
&=&(n-1)g^{jl}g^{pr}g^{qs}\nabla_j\nabla_l\left(\frac{\partial
g_{pq}}{\partial t}\frac{\partial g_{rs}}{\partial t}\right)\vspace{2mm}\\
\nonumber
&=&2(n-1)g^{jl}g^{pr}g^{qs}\left[\nabla_j\nabla_l(\frac{\partial
g_{pq}}{\partial t})\frac{\partial g_{rs}}{\partial
t}+\nabla_j(\frac{\partial g_{pq}}{\partial
t})\nabla_l(\frac{\partial g_{rs}}{\partial t})\right]\vspace{2mm}\\
\nonumber &=&2(n-1)g^{jl}g^{rs}\frac{\partial g_{ps}}{\partial
t}\dfrac{\partial }{\partial t}\frac{\partial
\Gamma_{lr}^p}{\partial x^j}+2(n-1)g^{jl}g^{rs}\frac{\partial
g_{pr}}{\partial t}\dfrac{\partial }{\partial t}\frac{\partial
\Gamma_{ls}^p}{\partial x^j}\vspace{2mm}\\ \nonumber
&&+2(n-1)g^{jl}g^{pq}g^{rs}\dfrac{\partial^2 g_{pr}}{\partial
t\partial x^j}\dfrac{\partial^2 g_{qs}}{\partial t\partial x^l}\vspace{2mm}\\
&=&4(n-1)g^{ik}g^{jl}\frac{\partial g_{ip}}{\partial
t}\dfrac{\partial }{\partial t}\frac{\partial
\Gamma_{kl}^p}{\partial
x^j}+2(n-1)g^{jl}g^{pq}g^{rs}\dfrac{\partial^2 g_{pr}}{\partial
t\partial x^j}\dfrac{\partial^2 g_{qs}}{\partial t\partial x^l}.
\end{eqnarray}
It follows from (6.7)-(6.14) that
\begin{eqnarray}
\nonumber\dfrac{\partial ^2R}{\partial t^2}&=& \Delta
R+2|Ric|^2\vspace{2mm}\\ \nonumber&&+\left[4\dfrac{\partial
g^{ik}}{\partial t}\dfrac{\partial R_{ik}}{\partial
t}-2\dfrac{\partial g^{ik}}{\partial t}\dfrac{\partial
g^{jl}}{\partial t}R_{ijkl}+4g^{ir}g^{ks}g^{pq}\dfrac{\partial
g_{pr}}{\partial
t}\dfrac{\partial g_{qs}}{\partial t}R_{ik}\right.\vspace{2mm}\\
\nonumber&&-2a\dfrac{\partial g^{ik}}{\partial t}\dfrac{\partial
R_{ik}}{\partial t}+a\dfrac{\partial g^{ik}}{\partial
t}\dfrac{\partial g^{jl}}{\partial
t}R_{ijkl}-2ag^{ir}g^{ks}g^{pq}\dfrac{\partial g_{pr}}{\partial
t}\dfrac{\partial g_{qs}}{\partial t}R_{ik}\vspace{2mm}\\
\nonumber&&\left.+(bg^{pq}\dfrac{\partial g_{pq}}{\partial
t}+d)\dfrac{\partial R}{\partial t}-\left(eg^{pq}\dfrac{\partial
g_{pq}}{\partial t}+f(g^{pq}\dfrac{\partial g_{pq}}{\partial
t})^2+h\dfrac{\partial g^{pq}}{\partial t}\dfrac{\partial
g_{pq}}{\partial t}\right)R\right]\vspace{2mm}\\
\nonumber&&+\left[(4(n-1)h-2a)g^{ik}g^{jl}\dfrac{\partial
g_{ip}}{\partial t}\dfrac{\partial }{\partial t}\frac{\partial
\Gamma_{kl}^p}{\partial x^j}-(2a+2b)\dfrac{\partial g^{ik}}{\partial
t}\dfrac{\partial }{\partial t}\frac{\partial
\Gamma_{kp}^p}{\partial x^i}\right.\vspace{2mm}\\
\nonumber&&\left.-\left(2bg^{pq}\dfrac{\partial g_{pq}}{\partial
t}+2(n-1)e+4(n-1)fg^{pq}\dfrac{\partial g_{pq}}{\partial
t}\right)g^{ik}\dfrac{\partial }{\partial t}\frac{\partial
\Gamma_{kp}^p}{\partial
x^i}\right]\vspace{2mm}\\
\nonumber&&+\left[2g^{ik}g^{jl}g_{pq}\dfrac{\partial
\Gamma_{ij}^p}{\partial t}\dfrac{\partial \Gamma_{kl}^q}{\partial
t}-2g^{ik}g^{jl}g_{pq}\dfrac{\partial \Gamma_{jl}^p}{\partial
t}\dfrac{\partial \Gamma_{ik}^q}{\partial t}\right.\vspace{2mm}\\
\nonumber&&\left.+ag^{ik}g^{jl}g^{pq}\left(\frac{\partial^2
g_{jp}}{\partial t\partial x^i}\frac{\partial^2 g_{kq}}{\partial
t\partial x^l}+\frac{\partial^2 g_{jp}}{\partial t\partial
x^l}\frac{\partial^2 g_{kq}}{\partial t\partial
x^i}-2\frac{\partial^2 g_{ip}}{\partial t\partial
x^j}\frac{\partial^2 g_{kq}}{\partial t\partial x^l}\right)\right.\vspace{2mm}\\
\nonumber&&\left.+bg^{ik}g^{jl}g^{pq}\left(\dfrac{\partial^2
g_{pq}}{\partial t\partial x^i}\dfrac{\partial^2 g_{jk}}{\partial
t\partial x^l}+\dfrac{\partial^2 g_{pq}}{\partial t\partial
x^l}\dfrac{\partial^2 g_{jk}}{\partial t\partial
x^i}-2\dfrac{\partial^2 g_{pq}}{\partial t\partial
x^j}\dfrac{\partial^2 g_{ik}}{\partial t\partial x^l}\right)\right.\vspace{2mm}\\
\nonumber&& \left.-2(n-1)fg^{jl}g^{pq}g^{rs}\dfrac{\partial^2
g_{pq}}{\partial t\partial x^j}\dfrac{\partial^2 g_{rs}}{\partial
t\partial x^l}+2(n-1)hg^{jl}g^{pq}g^{rs}\dfrac{\partial^2
g_{pr}}{\partial t\partial x^j}\dfrac{\partial^2 g_{qs}}{\partial
t\partial x^l}\right]\vspace{2mm}\\
\nonumber&=& \Delta R+2|Ric|^2\vspace{2mm}\\
\nonumber&&+\left[4\dfrac{\partial g^{ik}}{\partial
t}\dfrac{\partial R_{ik}}{\partial t}-2\dfrac{\partial
g^{ik}}{\partial t}\dfrac{\partial g^{jl}}{\partial
t}R_{ijkl}+4g^{ir}g^{ks}g^{pq}\dfrac{\partial g_{pr}}{\partial
t}\dfrac{\partial g_{qs}}{\partial t}R_{ik}\right.\vspace{2mm}\\
\nonumber&&-2a\dfrac{\partial g^{ik}}{\partial t}\dfrac{\partial
R_{ik}}{\partial t}+a\dfrac{\partial g^{ik}}{\partial
t}\dfrac{\partial g^{jl}}{\partial
t}R_{ijkl}-2ag^{ir}g^{ks}g^{pq}\dfrac{\partial g_{pr}}{\partial
t}\dfrac{\partial g_{qs}}{\partial t}R_{ik}\vspace{2mm}\\
\nonumber&&\left.+(bg^{pq}\dfrac{\partial g_{pq}}{\partial
t}+d)\dfrac{\partial R}{\partial t}-\left(eg^{pq}\dfrac{\partial
g_{pq}}{\partial t}+f(g^{pq}\dfrac{\partial g_{pq}}{\partial
t})^2+h\dfrac{\partial g^{pq}}{\partial t}\dfrac{\partial
g_{pq}}{\partial t}\right)R\right]\vspace{2mm}\\
\nonumber&&+\left[(4(n-1)h-2a)g^{ik}g^{jl}\dfrac{\partial
g_{ip}}{\partial t}\dfrac{\partial }{\partial t}\frac{\partial
\Gamma_{kl}^p}{\partial x^j}-(2a+2b)\dfrac{\partial g^{ik}}{\partial
t}\dfrac{\partial }{\partial t}\frac{\partial
\Gamma_{kp}^p}{\partial x^i}\right.\vspace{2mm}\\
\nonumber&&\left.-\left(2bg^{pq}\dfrac{\partial g_{pq}}{\partial
t}+2(n-1)e+4(n-1)fg^{pq}\dfrac{\partial g_{pq}}{\partial
t}\right)g^{ik}\dfrac{\partial }{\partial t}\frac{\partial
\Gamma_{kp}^p}{\partial
x^i}\right]\vspace{2mm}\\
\nonumber&&+\left(-2b-2(n-1)f-2\right)g^{ik}g^{jl}g^{pq}\dfrac{\partial^2
g_{ik}}{\partial t\partial x^p}\dfrac{\partial^2 g_{jl}}{\partial
t\partial x^q}\vspace{2mm}\\
\nonumber&&+(2b+8)g^{ik}g^{jl}g^{pq}\dfrac{\partial^2
g_{ik}}{\partial t\partial x^j}\dfrac{\partial^2 g_{pl}}{\partial
t\partial
x^q}+\left(2(n-1)h-2a+6\right)g^{ik}g^{jl}g^{pq}\dfrac{\partial^2
g_{ij}}{\partial t\partial x^p}\dfrac{\partial^2 g_{kl}}{\partial
t\partial x^q}\vspace{2mm}\\
&&+(a-8)g^{ik}g^{jl}g^{pq}\dfrac{\partial^2 g_{ij}}{\partial
t\partial x^l}\dfrac{\partial^2 g_{kp}}{\partial t\partial
x^q}+(a-4)g^{ik}g^{jl}g^{pq}\dfrac{\partial^2 g_{ip}}{\partial
t\partial x^j}\dfrac{\partial^2 g_{kl}}{\partial t\partial x^q}.
\end{eqnarray}

In (6.15), if we take
\begin{eqnarray}
a=2,\;b=-2,\;e=0,\;f=-\frac{2b}{4(n-1)}=\frac{1}{n-1},\;h=\frac{2a}{4(n-1)}=\frac{1}{n-1},
\end{eqnarray}
then we have
\begin{eqnarray}
\nonumber\dfrac{\partial ^2R}{\partial t^2}&=& \Delta
R+2|Ric|^2+(d-2g^{pq}\dfrac{\partial g_{pq}}{\partial
t})\dfrac{\partial R}{\partial
t}-\frac{1}{n-1}\left[(g^{pq}\dfrac{\partial g_{pq}}{\partial
t})^2+\dfrac{\partial g^{pq}}{\partial t}\dfrac{\partial
g_{pq}}{\partial t}\right]R\vspace{2mm}\\
\nonumber&&+\left[4g^{ik}g^{jl}g^{pq}\dfrac{\partial^2
g_{ik}}{\partial t\partial x^j}\dfrac{\partial^2 g_{pl}}{\partial
t\partial x^q}+4g^{ik}g^{jl}g^{pq}\dfrac{\partial^2 g_{ij}}{\partial
t\partial x^p}\dfrac{\partial^2 g_{kl}}{\partial
t\partial x^q}\right.\vspace{2mm}\\
&&\left.-6g^{ik}g^{jl}g^{pq}\dfrac{\partial^2 g_{ij}}{\partial
t\partial x^l}\dfrac{\partial^2 g_{kp}}{\partial t\partial
x^q}-2g^{ik}g^{jl}g^{pq}\dfrac{\partial^2 g_{ip}}{\partial t\partial
x^j}\dfrac{\partial^2 g_{kl}}{\partial t\partial x^q}\right].
\end{eqnarray}
In this case, the corresponding evolution equation reads
\begin{eqnarray}
\dfrac{\partial ^2g_{ij}}{\partial t^2}&=&-2R_{ij}+G_{ij},
\end{eqnarray}
where
\begin{equation}\begin{array}{lll}
G_{ij} & = &{\displaystyle 2g^{pq}\frac{\partial g_{ip}}{\partial
t}\frac{\partial g_{jq}}{\partial t}-2g^{pq}\frac{\partial
g_{pq}}{\partial t}\frac{\partial g_{ij}}{\partial
t}+d\frac{\partial g_{ij}}{\partial t}+}\vspace{2mm}\\ &
&{\displaystyle \frac{1}{n-1}\left[(g^{pq}\frac{\partial
g_{pq}}{\partial t})^2+(\frac{\partial g^{pq}}{\partial
t}\frac{\partial g_{pq}}{\partial t})\right]g_{ij}.}
\end{array}\end{equation}
Taking $d=-\tilde{d}$, we obtain from (6.18) and (6.19) that
\begin{equation}
\begin{array}{lll}
\dfrac{\partial^{2}g_{ij}}{\partial t^{2}}& = &
-2R_{ij}+2g^{pq}\dfrac{\partial g_{ip}}{\partial t}\dfrac{\partial
g_{jq}}{\partial t}-\left(\tilde{d}+2g^{pq}\dfrac{\partial
g_{pq}}{\partial
t}\right)\dfrac{\partial g_{ij}}{\partial t}+ \vspace{2mm}\\
& &\dfrac{1}{n-1}\left[\left(g^{pq}\dfrac{\partial g_{pq}}{\partial
t}\right)^2+ \dfrac{\partial g^{pq}}{\partial t}\dfrac{\partial
g_{pq}}{\partial t}\right]g_{ij},\end{array}
\end{equation}
where $\tilde{d}$ is a positive constant. Denoting $\tilde{d}$ by
$d$ in (6.20) gives the evolution equation (2.1) for the dissipative
hyperbolic geometric flow.

\begin{Theorem}
If we suppose that the evolution equation of the hyperbolic
geometric flow is defined by (6.18)-(6.19) on a manifold
$\mathscr{M}$, then the scalar curvature of the evolving metrics
satisfies the nonlinear wave equation (6.17) in the normal
coordinates.
\end{Theorem}
\begin{Remark}
If we take $a=b=d=e=f=h=0$, i.e., $G_{ij}\equiv 0$ in (6.1), then
(6.1) is nothing but the standard hyperbolic geometric flow (1.2)
(see \cite{kl}).
\end{Remark}
\begin{Remark}
For the evolution equation (6.17) of the scalar curvature, the last
term can be written in the covariant form as follow $$
g^{ik}g^{jl}g^{pq}\left(4\nabla_j\frac{\partial g_{ik}}{\partial
t}\nabla_q\frac{\partial g_{pl}}{\partial t}+4\nabla_p\frac{\partial
g_{ij}}{\partial t}\nabla_q\frac{\partial g_{kl}}{\partial
t}-6\nabla_l\frac{\partial g_{ij}}{\partial t}\nabla_q\frac{\partial
g_{kp}}{\partial t}-2\nabla_j\frac{\partial g_{ip}}{\partial
t}\nabla_q\frac{\partial g_{kl}}{\partial t}\right).
$$
\end{Remark}

\vskip 5mm

\noindent{\Large \textbf{Acknowledgements.}} The work of Kong was
supported in part by the NNSF of China (Grant No. 10671124) and the
NCET of China (Grant No. NCET-05-0390); the work of Liu was
supported in part by the NSF and NSF of China.

\end{document}